\documentclass[12pt]{article}
\usepackage{latexsym,amsfonts,amssymb}
\usepackage[dvips]{color}
\usepackage{colordvi,multicol}

\def \cal{\mathcal}

\newtheorem{thm}{Theorem}[section]
\newtheorem{cor}[thm]{Corollary}
\newtheorem{lem}[thm]{Lemma}

\date{}
\begin{document}
\title{\bf  The complement value problem for a class of second order elliptic integro-differential operators}
\author{}

\maketitle

\centerline{Wei Sun} \centerline{\small Department of Mathematics
and Statistics} \centerline{\small Concordia University} \centerline{\small
Montreal, H3G 1M8, Canada} \centerline{\small E-mail:
wei.sun@concordia.ca}

\vskip 1cm

%\vskip 1.4cm

\vskip 0.5cm \noindent{\bf Abstract}\quad  We consider the
complement value problem for a class of second order elliptic
integro-differential operators. Let $D$ be a bounded Lipschitz
domain of $\mathbb{R}^d$. Under mild conditions, we show that
there exists a unique bounded continuous weak solution to the
following equation
$$
\left\{\begin{array}{l}(\Delta+a^{\alpha}\Delta^{\alpha/2}+b\cdot\nabla+c+{\rm
div}
\hat{b})u+f=0\ \ {\rm in}\ D,\\
u=g\ \ {\rm on}\ D^c.
\end{array}\right.$$
 Moreover, we give an explicit
probabilistic representation of the solution. The recently
developed stochastic calculus for Markov processes associated with
semi-Dirichlet forms and heat kernel estimates play important
roles in our approach.

\smallskip

\noindent {\bf Keywords}\quad Complement value problem,
integro-differential operator, probabilistic representation,
semi-Dirichlet form, Fukushima type decomposition, heat kernel
estimate.

\smallskip

\noindent {AMS Subject Classification:} 35D30, 31C25, 60J75

\section{Introduction and main result}

Let $d\ge 1$ and $D$ be a bounded Lipschitz domain of
$\mathbb{R}^d$. Suppose that $0<\alpha<2$ and $p>d/2$. Let $a>0$,
$b=(b_1,\dots,b_d)^*$ satisfying $|b|\in L^{2p}(D;dx)$ if $d\ge 2$
and $|b|\in L^{\infty}(D;dx)$ if $d=1$, $c\in L^{p\vee 2}(D;dx)$,
$\hat{b}=(\hat{b}_1,\dots,\hat{b}_d)^*$ satisfying $|\hat{b}|\in
L^{2(p\vee 1)}(D;dx)$, $c+{\rm div} \hat{b}\le h$ for some $h\in
L^{p\vee 1}_+(D;dx)$ in the distribution sense, $f\in L^{4(p\vee
1)}(D;dx)$ and $g\in B_b(D^c)$. We consider the complement value
problem:
\begin{equation}\label{in1}
\left\{\begin{array}{l}(\Delta+a^{\alpha}\Delta^{\alpha/2}+b\cdot\nabla+c+{\rm
div}
\hat{b})u+f=0\ \ {\rm in}\ D,\\
u=g\ \ {\rm on}\ D^c.
\end{array}\right.
\end{equation}
The fractional Laplacian operator $\Delta^{\alpha/2}$ can be
written in the form
$$
\Delta^{\alpha/2}\phi(x)=\lim_{\varepsilon\rightarrow 0}{\cal
A}(d,-\alpha)\int_{\{|x-y|\ge
\varepsilon\}}\frac{\phi(y)-\phi(x)}{|x-y|^{d+\alpha}}dy,\ \
\phi\in C^{\infty}_c(\mathbb{R}^d),
$$
where ${\cal
A}(d,-\alpha):=\alpha2^{\alpha-1}\pi^{-d/2}\Gamma((d+\alpha)/2)\Gamma(1-\alpha/2)^{-1}$
and $C^{\infty}_c(\mathbb{R}^d)$ denotes the space of infinitely
differentiable functions on $\mathbb{R}^d$ with compact support.

The problem (\ref{in1}) is analogue of the Dirichlet problem for
second order elliptic integro-differential equations. For these
non-local equations, as opposed to the classical local case, the
function $g$ should be prescribed not only on the boundary
$\partial D$ but also in the whole complement $D^c$. The
complement value problem for non-local operators has many
applications, for example, in peridynamics \cite{AM,G,MD},
particle systems with long range interactions \cite{GL}, fluid
dynamics \cite{DG} and image processing \cite{GO}. This paper is
a continuation of our previous paper \cite{S}, which deals with
the case that $\hat{b}\equiv0$. We should point out that the lower
order term ${\rm div} \hat{b}$ in (\ref{in1}) is just a formal
writing since the vector field $\hat{b}$ is merely measurable
hence its divergence exists only in the distributional sense. Due
to the appearance of the ${\rm div} \hat{b}$, all the previous
known methods in solving the complement value problems such as
those in \cite{B} and \cite{S} ceased to work. In this paper, the
recently developed stochastic calculus for Markov processes
associated with semi-Dirichlet forms (see \cite{MSW, CMS}) will be
applied to overcome the difficulty caused by $\hat{b}$. In
particular, the Fukushima type decomposition for semi-Dirichlet
forms will play an important role in establishing the existence
and representation of the solution to the complement value
problem.

Denote $L:=\Delta+a^{\alpha}\Delta^{\alpha/2}+b\cdot\nabla$. By
setting $b = 0$ off $D$, we may assume that the operator $L$ is
defined on $\mathbb{R}^d$. By \cite[Theorem 1.4]{CH}, the
martingale problem for $(L, C^{\infty}_c(\mathbb{R}^d))$ is
well-posed for every initial value $x\in \mathbb{R}^d$. We use
$((X_t)_{t\ge 0},(P_x)_{x\in\mathbb{R}^d})$ to denote the strong
Markov process associated with $L$. Let $\rho>0$. Define
$$q_{\rho}(t,z)=t^{-d/2}\exp\left(-\frac{\rho|z|^2}{t}\right)+t^{-d/2}\wedge\frac{t}{|z|^{d+\alpha}},\
\ t>0,z\in\mathbb{R}^d.
$$ By
\cite[Theorems 1.2-1.4]{CH}, $X$ has a jointly continuous
transition density function $p(t,x,y)$ on $(0,\infty)\times
\mathbb{R}^d\times \mathbb{R}^d$, and for every $T>0$ there exist
positive constants $C_i,i=1,2,3,4$ such that
$$
C_1q_{C_2}(t,x-y)\le p(t,x,y)\le C_3 q_{C_4}(t,x-y),\ \ (t,x,y)\in
(0,T]\times\mathbb{R}^d\times\mathbb{R}^d.
$$

Define
$$
\left\{\begin{array}{l}{\cal
E}^0(\phi,\psi)=\int_{\mathbb{R}^d}\langle
\nabla\phi,\nabla\psi\rangle
dx+\frac{a^{\alpha}{\cal A}(d,-\alpha)}{2}\int_{\mathbb{R}^d}\int_{\mathbb{R}^d}\frac{(\phi(x)-\phi(y))(\psi(x)-\psi(y))}{|x-y|^{d+\alpha}}dxdy\\
\ \ \ \ \ \ \ \ \ \ \ \ \ \ \ -\int_{\mathbb{R}^d}\langle
b,\nabla \phi\rangle \psi dx,\ \ \phi,\psi\in D({\cal E}^0),\\
D({\cal E}^0)=W^{1,2}(\mathbb{R}^d)=\{ u\in L^2(\mathbb{R}^d;dx):|\nabla u|\in L^2(\mathbb{R}^d;dx)\}.
\end{array}
\right.
$$
By \cite[Lemma 2.1]{S}, we know that $({\cal E}^0,D({\cal E}^0))$
is a regular lower-bounded semi-Dirichlet form on
$L^2(\mathbb{R}^d;dx)$. Moreover, $((X_t)_{t\ge
0},(P_x)_{x\in\mathbb{R}^d})$ is the Hunt process associated with
$({\cal E}^0,D({\cal E}^0))$. By the assumption on $b$ and H\"older's inequality, we find that $|b|^2$ belongs to the Kato class. Then, we obtain by \cite[Chapter 7, Lemma 7.5]{Schechter} that there exists $\beta_0>0$ such that
$$
\int_{\mathbb{R}^d}|b|^2\phi^2dx\le\frac{1}{2}\int_{\mathbb{R}^d}|\nabla\phi|^2dx+\beta_0\int_{\mathbb{R}^d}|\phi|^2dx,\
\ \forall \phi\in W^{1,2}(\mathbb{R}^d).
$$
For $\phi\in C^{\infty}_c(\mathbb{R}^d)$, we have $L\phi\in
L^2(\mathbb{R}^d;dx)$ (cf. \cite[the proof of Lemm 2.1]{S}).
Define
$$\gamma:=\beta_0+1,
$$ and
$$ {\cal E}^0_{\gamma}(\phi,\psi)={\cal
E}^0(\phi,\psi)+\gamma(\phi,\psi),\ \ \phi,\psi\in D({\cal E}^0).
$$
Hereafter, $(\cdot,\cdot)$ denotes the inner product of
$L^2(\mathbb{R}^d;dx)$. Then, $({\cal E}^0_{\gamma},D({\cal
E}^0))$ is a regular semi-Dirichlet form on
$L^2(\mathbb{R}^d;dx)$.

Although $({\cal E}^0,D({\cal E}^0))$ is only a lower-bounded
semi-Dirichlet form, by replacing 1 and $({\cal E}^0_1,D({\cal
E}^0))$ with $\gamma$ and $({\cal E}^0_{\gamma},D({\cal E}^0))$
respectively in the proof of \cite[Theorem 1.4]{MSW}, we can
obtain the Fukushima type decomposition for $({\cal E}^0,D({\cal
E}^0))$. Let $\phi\in W^{1,2}(\mathbb{R}^d)$. We use
$\tilde{\phi}$ to denote an ${\cal E}^0$-quasi-continuous version
of $\phi$. Note that
$\mu_{\phi}(dx):=\int_{\mathbb{R}^d}\frac{(\phi(x)-\phi(y))^2}{|x-y|^{d+\alpha}}dydx$
is a finite measure on $\mathbb{R}^d$ and $X$ has infinite
lifetime. Then, $\phi$ admits a unique Fukushima type
decomposition
$$ \tilde\phi(X_t)-\tilde\phi(X_0)=M^{\phi}_t+N^{\phi}_t,\ \ t\ge
0,
$$
where $(M^{\phi}_t)_{t\ge 0}$ is a locally square integrable
martingale additive functional and $(N^{\phi}_t)_{t\ge 0}$ is a
continuous additive functional locally of zero quadratic
variation.

 By the Lax-Milgram theorem, for any vector
field $\eta\in L^2(\mathbb{R}^d;dx)$, there exists a unique
$\eta^H\in W^{1,2}(\mathbb{R}^d)$ such that
\begin{equation}\label{bbb}
\int_{\mathbb{R}^d}\langle \eta,\nabla\phi\rangle dx={\cal
E}^0_{\gamma}(\eta^H,\phi),\ \ \forall \phi\in
W^{1,2}(\mathbb{R}^d).
\end{equation}
Similar to \cite[Lemma 2.2]{CMS1}, we can show that
\begin{equation}\label{hatt}
\eta_n\rightarrow\eta\ {\rm in}\ L^2(\mathbb{R}^d;dx)\ {\rm as}\
n\rightarrow\infty\Longrightarrow\eta^H_n\rightarrow\eta^H\ {\rm
in}\ W^{1,2}(\mathbb{R}^d)\ {\rm as}\ n\rightarrow\infty.
\end{equation}
 Moreover, for $\phi\in C^{\infty}_c(\mathbb{R}^d)$, we have
\begin{equation}\label{smooth}
\int_0^t{\rm
div}\phi(X_s)ds=N^{\phi^H}_t-\gamma\int_0^t\phi^H(X_s)ds,\ \ t\ge
0. \end{equation}

Define
$$
e(t):=e^{\int_0^tc(X_s)ds+N^{\hat{b}^H}_t-\gamma\int_0^t\hat{b}^H(X_s)ds},\
\ t\ge0,
$$
and $\tau:=\inf\{t>0: X_t\in D^c\}$. Denote $W^{1,2}(D)=\{ u\in L^2(D;dx):|\nabla u|\in L^2(D;dx)\}$, $W^{1,2}_0(D)=\{ u\in W^{1,2}(D):\exists\{u_n\}_{n\in\mathbb{N}}\subset C_c^{\infty}(D)\ {\rm such\ that}\ u_n\rightarrow u\ {\rm in}\ W^{1,2}(D)\}$, and
$$
W^{1,2}_{loc}(D):=\{u: u\phi\in W^{1,2}_0(D)\ {\rm for\ any}\
\phi\in C^{\infty}_c(D)\}.
$$

The main result of this paper
is the following theorem.

\begin{thm}\label{thm1}
There exists $M>0$ such that if $\|h\|_{L^{p\vee 1}}\le M$, then
for any $f\in L^{4(p\vee 1)}(D;dx)$ and $g\in B_b(D^c)$, there
exists a unique $u\in B_b(\mathbb{R}^d)$ satisfying $u|_D\in
W^{1,2}_{loc}(D)\cap C(D)$ and
$$
\left\{\begin{array}{l}(\Delta+a^{\alpha}\Delta^{\alpha/2}+b\cdot\nabla+c+{\rm div} \hat{b})u+f=0\ \ {\rm in}\ D,\\
u=g\ \ {\rm on}\ D^c.
\end{array}\right.$$
Moreover, $u$ has the expression
$$
u(x)=E_x\left[e(\tau)g(X_{\tau})+\int_0^{\tau}e(s)f(X_s)ds\right]\
\ {\rm for\ q.e.}\ x\in D.
$$
In addition, if $g$ is continuous at $z\in \partial D$ then
$$
\lim_{x\rightarrow z}u(x)=u(z).
$$
\end{thm}

Hereafter
 $(\Delta+a^{\alpha}\Delta^{\alpha/2}+b\cdot\nabla+c+{\rm div} \hat{b})u+f=0$ is
understood in the distribution sense: for any $\phi\in
C^{\infty}_c(D)$,
\begin{eqnarray}\label{1}
& &\int_{D}\langle \nabla u,\nabla\phi\rangle
dx+\frac{a^{\alpha}{\cal A}(d,-\alpha)}{2}\int_{\mathbb{R}^d}\int_{\mathbb{R}^d}\frac{(u(x)-u(y))(\phi(x)-\phi(y))}{|x-y|^{d+\alpha}}dxdy\nonumber\\
& &\ \ \ \ \ \ \ \ \ \ -\int_{D}\langle b,\nabla u\rangle \phi
dx-\int_{D}cu\phi dx+\int_{D}\langle \hat{b},\nabla (u\phi)\rangle
dx-\int_{D}f\phi dx=0.
\end{eqnarray}
Note that the double integral appearing in (\ref{1}) is
well-defined for any $u\in B_b(\mathbb{R}^d)$ with $u|_D\in
W^{1,2}_{loc}(D)$ and $\phi\in C^{\infty}_c(D)$.

As a direct consequence of Theorem \ref{thm1}, we have the
following corollary.

\begin{cor}
If $c+{\rm div} \hat{b}\le 0$, then for any $f\in L^{4(p\vee
1)}(D;dx)$ and $g\in B_b(D^c)$ satisfying $g$ is continuous on
$\partial D$, there exists a unique $u\in B_b(\mathbb{R}^d)$ such
that $u$ is continuous on $\overline{D}$, $u|_D\in
W^{1,2}_{loc}(D)$, and
$$
\left\{\begin{array}{l}(\Delta+a^{\alpha}\Delta^{\alpha/2}+b\cdot\nabla+c+{\rm div} \hat{b})u+f=0\ \ {\rm in}\ D,\\
u=g\ \ {\rm on}\ D^c.
\end{array}\right.$$
Moreover, $u$ has the expression
$$
u(x)=E_x\left[e(\tau)g(X_{\tau})+\int_0^{\tau}e(s)f(X_s)ds\right]\
\ {\rm for\ q.e.}\ x\in D.
$$
\end{cor}

The proof of Theorem \ref{thm1} will be given in Section 3. In the
next section, we first present two lemmas, which will be used to
prove the continuity of the weak solution.

\section{Two lemmas}\setcounter{equation}{0}

 Throughout this paper, we denote by $C$ (or $C_i$) a generic fixed
strictly positive constant, whose value can change from line to
line. Define $B_r(x_0):=\{x\in\mathbb{R}^d:|x-x_0|<r\}$ for
$x_0\in \mathbb{R}^d$ and $r>0$.

\begin{lem}\label{Holder}
Suppose that $u\in B_b(\mathbb{R}^d)$ satisfying $u|_D\in
W^{1,2}_{loc}(D)$ and
$(\Delta+a^{\alpha}\Delta^{\alpha/2}+b\cdot\nabla+c+{\rm div}
\hat{b})u+f=0$ in $D$. Then, $u|_D$ has a locally H\"older
continuous version.
\end{lem}

\noindent {\bf Proof.}\ \ It suffices to assume that $d\ge 2$. Set
$\overline{u}=u+2\|u\|_{\infty}$,
$\overline{b}=b-2\|u\|_{\infty}\overline{u}^{-1}\hat{b}$,
$\overline{c}=(1-2\|u\|_{\infty}\overline{u}^{-1})c+f\overline{u}^{-1}$
and
$\overline{\hat{b}}=(1-2\|u\|_{\infty}\overline{u}^{-1})\hat{b}$.
Then,
$$
(\Delta+a^{\alpha}\Delta^{\alpha/2}+\overline{b}\cdot\nabla+\overline{c}+{\rm
div} \overline{\hat{b}})\overline{u}=0\ \ {\rm in}\ D.
$$
Thus, to prove Lemma \ref{Holder}, we may assume without loss of
generality that $f\equiv0$.

The  proof given below follows Kassmann \cite{Kass}. Denote ${\cal
L}:=\Delta+a^{\alpha}\Delta^{\alpha/2}+b\cdot\nabla+c+{\rm div}
\hat{b}$.

\noindent (i)\ \ Assume that $0<\rho<r<1$, $x_0\in\mathbb{R}^d$
and $\omega\in L^{q/2}_{loc}(B_{2r}(x_0))$ for some $q>d$. Suppose
that $v\in W^{1,2}(B_{2r}(x_0))$ is nonnegative in $\mathbb{R}^d$
and satisfies $(-{\cal L}v,\phi)\ge(\omega,\phi)$ for any
nonnegative $\phi\in C_c^{\infty}(B_{2r}(x_0))$ and $v(x)\ge
\varepsilon$ for almost all $x\in B_{2r}(x_0)$ and some
$\varepsilon>0$. Let $\theta\in C^{\infty}_c(B_{5r/4}(x_0))$ be a
function with ${\rm supp}(\theta)=\overline{B_{r+\rho}(x_0)}$,
$\theta(x)=1$ for $x\in B_{r}(x_0)$ and $\theta(x)>0$ for $x\in
B_{r+\rho}(x_0)$, $\|\nabla\theta\|_{\infty}\le \iota\rho^{-1}$ for
some constant $\iota>0$ and $\|\theta\|_{\infty}\le 1$.

We have
\begin{eqnarray*}
(\omega,-\theta^2v^{-1})&\ge&((\Delta+b\cdot\nabla+c+{\rm div}
\hat{b})v,\theta^2v^{-1})\\
& &+\frac{a^{\alpha}{\cal
A}(d,-\alpha)}{2}\int_{\mathbb{R}^d}\int_{\mathbb{R}^d}\frac{(v(y)-v(x))(\theta^2v^{-1}(x)-\theta^2v^{-1}(y))}{|x-y|^{d+\alpha}}dxdy.
\end{eqnarray*}
By \cite[(5.10)]{TR} and \cite[(5.7)--(5.9)]{Kass}, there exists a
constant $C>0$ which is independent of $v,x_0,r,\rho,\omega$ and
$\varepsilon$ such that
\begin{equation}\label{JN}
\int_{\mathbb{R}^d}\theta^2|\nabla(\log v)|^2dx\le
C\rho^{-2}|B_{r+\rho}(x_0)|+\varepsilon^{-1}\|\omega\|_{L^{q/2}(B_{r+\rho}(x_0))}\|1\|_{L^{q/q-2}(B_{r+\rho}(x_0))}.
\end{equation}

\noindent (ii) Assume that $0<R<1$, $x_0\in\mathbb{R}^d$ and
$\omega\in L^{q/2}_{loc}(B_{5R/4}(x_0))$ for some $q>d$. Suppose
that $v\in W^{1,2}(B_{5R/4}(x_0))$ is nonnegative in
$\mathbb{R}^d$ and satisfies $(-{\cal L}v,\phi)\ge(\omega,\phi)$
for any nonnegative $\phi\in C_c^{\infty}(B_{5R/4}(x_0))$ and
$v(x)\ge \varepsilon$ for almost all $x\in B_{5R/4}(x_0)$ and some
$\varepsilon>\frac{1}{4}R^{\frac{2(q-d)}{q}}\|\omega\|_{L^{q/2}(B_{9R/8}(x_0))}$.
By (\ref{JN}), similar to \cite[Lemma 5.11]{Kass}, we can show
that there exist $\bar{p}\in (0,1)$ and $\iota>0$ such that
\begin{equation}\label{ku1}
\left(|B_R(x_0)|^{-1}\int_{B_R(x_0)}v^{{\bar{p}}}dx\right)^{1/{\bar{p}}}\le
\iota\left(|B_R(x_0)|^{-1}\int_{B_R(x_0)}v^{-{\bar{p}}}dx\right)^{-1/{\bar{p}}},
\end{equation}
where $\iota$ and $\bar{p}$ are independent of $v,x_0,R$ and
$\varepsilon$.

\noindent (iii)\ \ Assume that $0<8\rho<R<1-\rho$,
$x_0\in\mathbb{R}^d$ and $\omega\in L^{q/2}_{loc}(B_{5R/4}(x_0))$
for some $q>d$. Suppose that $v\in W^{1,2}(B_{5R/4}(x_0))$
satisfying $(-{\cal L}v,\phi)\ge(\omega,\phi)$ for any nonnegative
$\phi\in C_c^{\infty}(B_{R}(x_0))$ and $v(x)\ge \varepsilon$ for
almost all $x\in B_{R}(x_0)$ and some
$\varepsilon>R^{\frac{2(q-d)}{q}}\|\omega\|_{L^{q/2}(B_{9R/8}(x_0))}$.
Let $\theta\in C^{\infty}_c(B_{9R/8}(x_0))$ be a function with
${\rm supp}(\theta)=\overline{B_{R+\rho}(x_0)}$, $\theta(x)=1$ for
$x\in B_{R}(x_0)$ and $\theta(x)>0$ for $x\in B_{R+\rho}(x_0)$,
$\|\nabla\theta\|_{\infty}\le \iota\rho^{-1}$ for some constant $\iota>0$
and $\|\theta\|_{\infty}\le 1$.

Let $\varrho>1$. We have
\begin{eqnarray*}
(\omega,-\theta^{\varrho+1}v^{-\varrho})&\ge&((\Delta+b\cdot\nabla+c+{\rm
div}
\hat{b})v,\theta^{\varrho+1}v^{-\varrho})\\
& &+\frac{a^{\alpha}{\cal
A}(d,-\alpha)}{2}\int_{\mathbb{R}^d}\int_{\mathbb{R}^d}\frac{(v(y)-v(x))(\theta^{\varrho+1}v^{-\varrho}(x)-\theta^{\varrho+1}v^{-\varrho}(y))}{|x-y|^{d+\alpha}}dxdy.
\end{eqnarray*}
By \cite[(5.10)]{TR} and \cite[(5.14) and (5.15)]{Kass}, there
exists a constant $C>0$ which is independent of
$v,x_0,R,\rho,\varrho$ and $\varepsilon$ such that
\begin{equation}\label{ha}
\|v^{-1}\|^{\varrho-1}_{L^{(\varrho-1)\frac{d^*}{d^*-2}}(B_R(x_0))}\le
C(\max\{\varrho-1,(\varrho-1)^2\})\rho^{-2}\|v^{-1}\|^{\varrho-1}_{L^{\varrho-1}(B_{R+\rho}(x_0)},
\end{equation}
where $d^*=d$ if $d\ge 3$ and $2<d^*<q$ if $d=2$.

\noindent (iv) Assume that  $0<R<1/2$, $0<\mu<1<\Theta$,
$x_0\in\mathbb{R}^d$ and $\omega\in L^{q/2}_{loc}(B_{\Theta
R}(x_0))$ for some $q>d$. Suppose that $v\in W^{1,2}(B_{\Theta
R}(x_0))$ satisfying $(-{\cal L}v,\phi)\ge(\omega,\phi)$ for any
nonnegative $\phi\in C_c^{\infty}(B_{\Theta R}(x_0))$ and $v(x)\ge
\varepsilon$ for almost all $x\in B_{\Theta R}(x_0)$ and some
$\varepsilon>(\Theta
R)^{\frac{2(q-d)}{q}}\|\omega\|_{L^{q/2}(B_{\frac{1+3\Theta}{4}R}(x_0))}$.
By (\ref{ha}), similar to \cite[Corollary 5.13]{Kass}, we can show
that for any $\varrho_0>0$,
\begin{equation}\label{ku2}
\inf_{x\in B_{\mu R}(x_0)}v(x)\ge
\iota\left(|B_R(x_0)|^{-1}\int_{B_R(x_0)}v^{-\varrho_0}dx\right)^{-1/\varrho_0},
\end{equation}
where $\iota>0$ is independent of $v,x_0,R$ and $\varepsilon$.

\noindent (v)\ \ By (\ref{ku1}) and (\ref{ku2}), similar to
\cite[Corollary 5.9]{Kass}, we can establish the following weak
Harnack inequality:

There exist positive constants $\iota_1,\iota_2$ and $\varrho_0$ such that
for any $x_0\in\mathbb{R}^d$, $R\in (0,1)$, $v\in
W^{1,2}(B_{R}(x_0))$ satisfying $v\ge 0$ in $B_R(x_0)$ and
$(-{\cal L}v,\phi)\ge0$ for any nonnegative $\phi\in
C_c^{\infty}(B_{R}(x_0))$, we have
$$
\inf_{B_{R/4}(x_0)}v\ge
\iota_1\left(|B_{R/2}(x_0)|^{-1}\int_{B_{R/2}(x_0)}v^{\varrho_0}dx\right)^{1/\varrho_0}-\iota_2R^2\sup_{x\in
B_{R/2}(x_0)}\int_{\mathbb{R}^d\backslash
B_R(x_0)}\frac{v^{-}(z)}{|x-z|^{d+\alpha}}dz,
$$
where $v^{-}(z)=-(v(z)\wedge 0)$. The proof is therefore complete
by the weak Harnack inequality and \cite[Corollary 4.2]{Kass}.
\hfill\fbox

\vskip 0.2cm Suppose
 $d\ge 2$. We define $q=\frac{p}{p-1}$. Then $\frac{1}{p}+\frac{1}{q}=1$ and $1<q<\frac{d}{d-2}$. We choose $\beta$
 such that
$$\label{new42}\frac{d}{2}-1<\beta<\frac{d}{2q}.$$ Let
$M_1>0$ be a constant satisfying
\begin{eqnarray}\label{June11}
e^{|x|}\ge M_1|x|^{\beta},\ \ \forall x\in\mathbb{R}^d.
\end{eqnarray}
We define $q^*=\frac{2p}{2p-1}$. Then
$\frac{1}{2p}+\frac{1}{q^*}=1$ and $1<q^*<\frac{d}{d-1}$. We
choose $0<\gamma<1$ such that
$$
q^*<\frac{d}{d-\gamma}.
$$
Let $M_2>0$ be a constant satisfying
\begin{eqnarray}\label{June1}
e^{|x|}\ge M_2|x|^{(d-\gamma)/2},\ \ \forall x\in\mathbb{R}^d.
\end{eqnarray}
We choose $0<\delta<(2-\alpha)/2$. Let $M_3>0$ be a constant
satisfying
\begin{eqnarray}\label{July}
e^{|x|}\ge M_3|x|^{1/6},\ \ \forall x\in\mathbb{R}^d.
\end{eqnarray}
Define $\varsigma=\sup_{x\in D}|x|$.

\begin{lem}\label{lemma23} Let $C,R$ be two positive constants and $\mu$ be a function on $\mathbb{R}^d$ with ${\rm supp}[\mu]\subset B_R(0)$.

(i) Suppose $d\ge 2$. Then, there exist positive constants
${\cal C}_1,{\cal C}_2$ which are independent of $\mu$ such that for any $t>0$
and $x\in D$,
\begin{eqnarray*}&&\int_0^t\int_{y\in
\mathbb{R}^d}\left(s^{-d/2}\exp\left(-\frac{C|x-y|^2}{s}\right)+s^{-d/2}\wedge\frac{s}{|x-y|^{d+\alpha}}\right)|\mu
(y)|dyds\\
&\le& {\cal C}_1(t^{\beta+1-d/2}+t^{\delta})\left(\int_{y\in
\mathbb{R}^d}|\mu(y)|^{p}dy\right)^{1/p}\end{eqnarray*} and
\begin{eqnarray*}&&\int_0^t\int_{y\in
\mathbb{R}^d}\left(s^{-(d+1)/2}\exp\left(-\frac{C|x-y|^2}{s}\right)+s^{-(d+1)/2}\wedge\frac{s}{|x-y|^{d+1+\alpha}}\right)|\mu
(y)|dyds\\
&\le &{\cal C}_2(t^{(1-\gamma)/2}+t^{\delta})\left(\int_{y\in
\mathbb{R}^d}|\mu(y)|^{2p}dy\right)^{1/(2p)}.
\end{eqnarray*}

(ii)  Suppose $d=1$. Then, for any $t>0$ and $x\in D$,
\begin{eqnarray}\label{Jane}
&&\int_0^t\int_{-\infty}^{\infty}\left(s^{-1/2}\exp\left(-\frac{C|x-y|^2}{s}\right)+s^{-1/2}\wedge\frac{s}{|x-y|^{1+\alpha}}\right)|\mu
(y)|dyds\nonumber\\
&\le&4t^{1/2}\int_{-\infty}^{\infty}|\mu (y)|dy,
\end{eqnarray}
and there exists a positive constant ${\cal C}_3$ which is independent of
$\mu$ such that for any $t>0$ and $x\in D$,
\begin{eqnarray*}&&\int_0^t\int_{-\infty}^{\infty}\left(s^{-1}\exp\left(-\frac{C|x-y|^2}{s}\right)+s^{-1}\wedge\frac{s}{|x-y|^{2+\alpha}}\right)|\mu
(y)|dyds\\
&\le&{\cal C}_3(t^{1/6}+t^{\delta})\left(\int_{-\infty}^{\infty}|\mu(y)|^{2}dy\right)^{1/2}.
\end{eqnarray*}
\end{lem}

\noindent {\bf Proof.} We first consider the case that $d\ge 2$. We have \begin{equation}\label{zxzx}
s^{-d/2}\wedge\frac{s}{|x-y|^{d+\alpha}}\le s^{-d/2}\le
e^{C}s^{-d/2}\exp\left(-\frac{C|x-y|^2}{s}\right)\ \ {\rm if}\
|x-y|^2<s.
\end{equation}
By (\ref{June11}) and (\ref{zxzx}), we get
\begin{eqnarray*}
&&\int_0^t\int_{y\in
\mathbb{R}^d}\left(s^{-d/2}\exp\left(-\frac{C|x-y|^2}{s}\right)+s^{-d/2}\wedge\frac{s}{|x-y|^{d+\alpha}}\right)|\mu
(y)|dyds\\
&\le&\int_0^t\int_{y\in
\mathbb{R}^d}(1+e^C)s^{-d/2}\exp\left(-\frac{C|x-y|^2}{s}\right)|\mu
(y)|dyds\\
&&+\int_0^t\int_{\{|x-y|\ge\sqrt{s}\}}\frac{s}{|x-y|^{d+\alpha}}|\mu
(y)|dyds\\
&\le&\int_0^t\int_{y\in
\mathbb{R}^d}\frac{1+e^C}{M_1s^{d/2}(C|x-y|^2/s)^{\beta}}|\mu(y)|dyds+t^{\delta}\int_{{y\in
\mathbb{R}^d}}\int_0^{|x-y|^2}s^{1-\delta}ds\frac{1}{|x-y|^{d+\alpha}}|\mu
(y)|dy\nonumber\\
&=&\int_0^t\frac{1+e^C}{C^\beta M_1s^{d/2-\beta}}ds\int_{y\in
\mathbb{R}^d}\frac{|\mu(y)|}{|x-y|^{2\beta}}dy+t^{\delta}\int_{{y\in
\mathbb{R}^d}}\frac{1}{(2-\delta)|x-y|^{d+\alpha-4+2\delta}}|\mu
(y)|dy\nonumber\\
&\le&\int_0^t\frac{1+e^C}{C^\beta
M_1s^{d/2-\beta}}ds\left(\int_{y\in B_R(0)}\frac{1}{|x-y|^{2\beta
q}}dy\right)^{1/q}\left(\int_{y\in \mathbb{R}^d}|\mu(y)|^{p}
dy\right)^{1/p}
\nonumber\\
&&+\frac{t^{\delta}}{2-\delta}\left(\int_{y\in
B_R(0)}\frac{1}{|x-y|^{(d+\alpha-4+2\delta)q}}dy\right)^{1/q}\left(\int_{y\in
\mathbb{R}^d}|\mu(y)|^{p} dy\right)^{1/p}
\nonumber\\
 &=&\frac{(1+e^C)t^{\beta+1-d/2}}{C^\beta
M_1(\beta+1-d/2)}\left(C_1\int_0^{R+\varsigma} r^{d-2\beta
q-1}dr\right)^{1/q}\left(\int_{y\in
\mathbb{R}^d}|\mu(y)|^{p}dy\right)^{1/p}\nonumber\\
&&+ \frac{t^{\delta}}{2-\delta}\left(C_2\int_0^{R+\varsigma}
r^{d-(d+\alpha-4+2\delta)q-1}dr\right)^{1/q}\left(\int_{y\in
\mathbb{R}^d}|\mu(y)|^{p}dy\right)^{1/p}\\
&:=&{\cal C}_1(t^{\beta+1-d/2}+t^{\delta})\left(\int_{y\in
\mathbb{R}^d}|\mu(y)|^{p}dy\right)^{1/p}.
\end{eqnarray*}
By (\ref{June1}) and (\ref{zxzx}), we get
\begin{eqnarray*}
&&\int_0^t\int_{y\in
\mathbb{R}^d}\left(s^{-(d+1)/2}\exp\left(-\frac{C|x-y|^2}{s}\right)+s^{-(d+1)/2}\wedge\frac{s}{|x-y|^{d+1+\alpha}}\right)|\mu
(y)|dyds\\
&\le&\int_0^t\int_{y\in
\mathbb{R}^d}(1+e^C)s^{-(d+1)/2}\exp\left(-\frac{C|x-y|^2}{s}\right)|\mu
(y)|dyds\\
&&+\int_0^t\int_{\{|x-y|\ge\sqrt{s}\}}\frac{s}{|x-y|^{d+1+\alpha}}|\mu
(y)|dyds\\
&\le&\int_0^t\int_{y\in
\mathbb{R}^d}\frac{1+e^C}{M_2s^{(d+1)/2}(C|x-y|^2/s)^{(d-\gamma)/2}}|\mu(y)|dyds+t^{\delta}\int_{y\in
\mathbb{R}^d}\int_0^{|x-y|^2}s^{1-\delta}ds\frac{1}{|x-y|^{d+1+\alpha}}|\mu
(y)|dy\nonumber\\
&=&\int_0^t\frac{1+e^C}{C^{(d-\gamma)/2}
M_2s^{(1+\gamma)/2}}ds\int_{y\in
\mathbb{R}^d}\frac{|\mu(y)|}{|x-y|^{d-\gamma}}dy+t^{\delta}\int_{y\in
\mathbb{R}^d}\frac{1}{(2-\delta)|x-y|^{d+\alpha-3+2\delta}}|\mu
(y)|dy\nonumber\\
&\le&\int_0^t\frac{1+e^C}{C^{(d-\gamma)/2}
M_2s^{(1+\gamma)/2}}ds\left(\int_{y\in
B_R(0)}\frac{1}{|x-y|^{(d-\gamma)q^*}}dy\right)^{1/q^*}\left(\int_{y\in
\mathbb{R}^d}|\mu(y)|^{2p} dy\right)^{1/(2p)}
\nonumber\\
&&+\frac{t^{\delta}}{2-\delta}\left(\int_{y\in
B_R(0)}\frac{1}{|x-y|^{(d+\alpha-3+2\delta)q^*}}dy\right)^{1/q^*}\left(\int_{y\in
\mathbb{R}^d}|\mu(y)|^{2p} dy\right)^{1/(2p)}
\nonumber\\
 &=&\frac{2(1+e^C)t^{(1-\gamma)/2}}{C^{(d-\gamma)/2}
M_2(1-\gamma)}\left(C_3\int_0^{R+\varsigma}
r^{d-(d-\gamma)q^*-1}dr\right)^{1/q^*}\left(\int_{y\in
\mathbb{R}^d}|\mu(y)|^{2p}dy\right)^{1/(2p)}\nonumber\\
&&+ \frac{t^{\delta}}{2-\delta}\left(\int_0^{R+\varsigma}
r^{d-(d+\alpha-3+2\delta)q^*-1}dr\right)^{1/q}\left(C_4\int_{y\in
\mathbb{R}^d}|\mu(y)|^{2p}dy\right)^{1/(2p)}\\
&:=&{\cal C}_2(t^{(1-\gamma)/2}+t^{\delta})\left(\int_{y\in
\mathbb{R}^d}|\mu(y)|^{2p}dy\right)^{1/(2p)}.
\end{eqnarray*}

We now consider the case that $d=1$. It is easy to see that
(\ref{Jane}) holds. By (\ref{July}) and (\ref{zxzx}), we get
\begin{eqnarray*}
&&\int_0^t\int_{-\infty}^{\infty}\left(s^{-1}\exp\left(-\frac{C|x-y|^2}{s}\right)+s^{-1}\wedge\frac{s}{|x-y|^{2+\alpha}}\right)|\mu
(y)|dyds\\
&\le&\int_0^t\int_{-\infty}^{\infty}(1+e^C)s^{-1}\exp\left(-\frac{C|x-y|^2}{s}\right)|\mu
(y)|dyds\\
&&+\int_0^t\int_{\{|x-y|\ge\sqrt{s}\}}\frac{s}{|x-y|^{2+\alpha}}|\mu
(y)|dyds\\
&\le&\int_0^t\int_{-\infty}^{\infty}\frac{1+e^C}{M_3s(C|x-y|^2/s)^{1/6}}|\mu(y)|dyds+t^{\delta}\int_{-\infty}^{\infty}\int_0^{|x-y|^2}s^{1-\delta}ds\frac{1}{|x-y|^{2+\alpha}}|\mu
(y)|dy\nonumber\\
&=&\int_0^t\frac{1+e^C}{C^{1/6}
M_3s^{5/6}}ds\int_{-\infty}^{\infty}\frac{|\mu(y)|}{|x-y|^{1/3}}dy+t^{\delta}\int_{-\infty}^{\infty}\frac{1}{(2-\delta)|x-y|^{\alpha-2+2\delta}}|\mu
(y)|dy\nonumber\\
&\le&\frac{6(1+e^C)t^{1/6}}{C^{1/6}
M_3}\left(\int_{-R}^R\frac{1}{|x-y|^{2/3}}dy\right)^{1/2}\left(\int_{-\infty}^{\infty}|\mu(y)|^2
dy\right)^{1/2}
\nonumber\\
&&+\frac{t^{\delta}}{2-\delta}\left(\int_{y\in
B_R(0)}\frac{1}{|x-y|^{2(\alpha-2+2\delta)}}dy\right)^{1/2}\left(\int_{-\infty}^{\infty}|\mu(y)|^{2}
dy\right)^{1/2}
\nonumber\\
 &=&\frac{6(1+e^C)t^{1/6}}{C^{1/6}
M_3}\left(C_5\int_0^{R+\varsigma}
r^{-2/3}dr\right)^{1/2}\left(\int_{-\infty}^{\infty}|\mu(y)|^{2}dy\right)^{1/2}\nonumber\\
&&+ \frac{t^{\delta}}{2-\delta}\left(C_6\int_0^{R+\varsigma}
r^{2(2-\alpha-2\delta)}dr\right)^{1/q}\left(\int_{-\infty}^{\infty}|\mu(y)|^{2}dy\right)^{1/2}\\
&:=&{\cal C}_3(t^{1/6}+t^{\delta})\left(\int_{-\infty}^{\infty}|\mu(y)|^{2}dy\right)^{1/2}.
\end{eqnarray*}
\hfill\fbox

\section{Proof of Theorem \ref{thm1}}\setcounter{equation}{0}

\subsection{Existence of weak solution}\label{dfgh}

Define
\begin{equation}\label{weak2}
u^*(x)=E_x\left[e(\tau)g(X_{\tau})+\int_0^{\tau}e(s)f(X_s)ds\right]\
\ {\rm for}\ x\in \mathbb{R}^d. \end{equation} By \cite[Lemma
2.5]{S}, there exists $C>0$  such that
\begin{equation}\label{12}
\sup_{x\in D}E_x\left[\int_0^{\tau}v(X_s)ds\right]\le
C\|v\|_{L^{p\vee 1}},\ \ \forall v\in L_+^{p\vee 1}(D).
\end{equation}
By Khasminskii's inequality and (\ref{12}), we find that there
exist positive constants $M$ and $\Upsilon$ such that for any
$v\in L^{p\vee 1}_+(D)$ satisfying $\|v\|_{L^{p\vee 1}}\le M$, we
have
\begin{equation}\label{14}
\sup_{x\in
D}E_x\left[e^{\int_0^{\tau}8v(X_s)ds}\right]\le\Upsilon.
\end{equation}
In particular, this implies that there exists $\nu>0$ such that
\begin{equation}\label{15}
\sup_{x\in D}E_x\left[e^{\nu\tau}\right]<\infty.
\end{equation}

Define
$$
J(x)=\frac{1_{\{|x|<1\}}e^{-\frac{1}{1-|x|^2}}}{\int_{\{|y|<1\}}e^{-\frac{1}{1-|y|^2}}dy},\
\ x\in \mathbb{R}^d.
$$
For $k\in \mathbb{N}$ and $x\in \mathbb{R}^d$, set
\begin{eqnarray*}
J_k(x)&=&k^dJ(kx),\\
\hat{b}_k(x)&=&\int_{\mathbb{R}^d}\hat{b}(x-y)J_k(y)dy,\\
c_k(x)&=&\int_{\mathbb{R}^d}c(x-y)J_k(y)dy,\\
h_k(x)&=&\int_{\mathbb{R}^d}h(x-y)J_k(y)dy.
\end{eqnarray*}
We have
\begin{equation}\label{may1}
\hat{b}_k\rightarrow\hat{b}\ \ {\rm in}\ \ L^2(\mathbb{R}^d;dx)\ \
{\rm as}\ \ k\rightarrow\infty,
\end{equation}
and
\begin{equation}\label{may2}
c_k\rightarrow c\ \ {\rm in}\ \ L^2(\mathbb{R}^d;dx)\ \ {\rm as}\
\ k\rightarrow\infty.
\end{equation}
Suppose that $\|h\|_{L^{p\vee 1}}\le M$. Since $c+{\rm
div}\hat{b}\le h$ implies that $c_k+{\rm div}\hat{b}_k\le h_k$ for
$k\in \mathbb{N}$, we obtain by (\ref{14}) that
\begin{equation}\label{may3}
\sup_{k\in\mathbb{N}}\sup_{x\in
D}E_x\left[e^{\int_0^{\tau}8(c_k+{\rm
div}\hat{b}_k)(X_s)ds}\right]\le\sup_{k\in\mathbb{N}}\sup_{x\in
D}E_x\left[e^{\int_0^{\tau}8h_k(X_s)ds}\right]\le\Upsilon.
\end{equation}

By (\ref{hatt}) and (\ref{may1}), we get
\begin{equation}\label{may4}
\hat{b}^H_k\rightarrow\hat{b}^H\ \ {\rm in}\ \
W^{1,2}(\mathbb{R}^d)\ \ {\rm as}\ \ k\rightarrow\infty.
\end{equation}
Further, by \cite[Theorem 1.15]{MSW} and \cite[Lemma A.6, Theorem
A.8 and Lemma A.9]{MMS}, similar to \cite[Corollary 5.2.1]{Fuku}
we can show that there exists a subsequence $\{k_l\}$ such that
for q.e. $x\in\mathbb{R}^d$,
\begin{equation}\label{may5}
P_x\left\{\lim_{l\rightarrow\infty}N^{\hat{b}^H_{k_l}}_t=N^{\hat{b}^H}_t\
{\rm uniformly\ on\ any\ finite\ interval\ of}\ t\right\}=1.
\end{equation}
For simplicity, we still use $\{k\}$ to denote the subsequence
$\{k_l\}$. By (\ref{smooth}), (\ref{12}),
(\ref{may2})--(\ref{may5}) and Fatou's lemma, we find that there
exists an ${\cal E}^0$-exceptional set $F\subset D$ such that for
$x\in D\backslash F$,
\begin{eqnarray}\label{div1}
E_x[e^2(\tau)]&\le&\liminf_{k\rightarrow\infty}E_x\left[e^{2\left[\int_0^{\tau}c_k(X_s)ds+N^{\hat{b}_k^H}_{\tau}-\gamma\int_0^{\tau}\hat{b}_k^H(X_s)ds\right]}\right]\nonumber\\
&=&\liminf_{k\rightarrow\infty}E_x\left[e^{\int_0^{\tau}2(c_k+{\rm div}\hat{b}_k)(X_s)ds}\right]\nonumber\\
&\le&\liminf_{k\rightarrow\infty}E_x\left[e^{\int_0^{\tau}2(h_k)(X_s)ds}\right]\nonumber\\
&\le&\Upsilon^{1/4},
\end{eqnarray}
and
\begin{eqnarray}\label{13}
& &E_x\left[\left(\int_0^{\tau}e(s)f(X_s)ds\right)^2\right]\nonumber\\
&\le&\liminf_{k\rightarrow\infty}E_x\left[\left(\int_0^{\tau}e^{\int_0^s(c_k+{\rm div}\hat{b}_k)(X_t)dt}f(X_s)ds\right)^2\right]\nonumber\\
&\le&\liminf_{k\rightarrow\infty}\left(E_x\left[\left(\int_0^{\tau}(\tau\vee 1) e^{\int_0^s2(c_k+{\rm div}\hat{b}_k)(X_t)dt}ds\right)^2\right]\right)^{1/2}\left(E_x\left[\left(\int_0^{\tau}\frac{1}{\tau\vee 1}f^2(X_s)ds\right)^2\right]\right)^{1/2}\nonumber\\
&\le&\liminf_{k\rightarrow\infty}\left(E_x\left[e^{\int_0^{\tau}4h_k(X_s)ds}\cdot\tau^2(\tau\vee 1)^2\right]\right)^{1/2}\left(E_x\left[\int_0^{\tau}f^4(X_s)ds\right]\right)^{1/2}\nonumber\\
&\le&C^{1/2}\liminf_{k\rightarrow\infty}\left(E_x\left[e^{\int_0^{\tau}8h_k(X_s)ds}\right]\right)^{1/4}\left(E_x\left[\tau^4(\tau\vee
1)^4\right]\right)^{1/4}\|f^4\|^{1/2}_{L^{p\vee
1}}\nonumber\\
&\le&C^{1/2}\Upsilon^{1/4}\left(E_x\left[\tau^4(\tau\vee
1)^4\right]\right)^{1/4}\|f^4\|^{1/2}_{L^{p\vee 1}}.
\end{eqnarray}
By (\ref{15}), (\ref{div1}) and (\ref{13}), we know that there
exists $M>0$ such that if $\|h\|_{L^{p\vee 1}}\le M$, then for any
$f\in L^{4(p\vee 1)}(D;dx)$ and $g\in B_b(D^c)$, $u^*$ is bounded
on $\mathbb{R}^d$ except for an ${\cal E}^0$-exceptional set.

For $k\in \mathbb{N}$, define
$$
e_k(t):=e^{\int_0^t(c_k+{\rm div}\hat b_k)(X_s)ds},\ \ t\ge0.
$$
By \cite[Theorem 1.1]{S}, we know that the unique bounded
continuous weak solution to the complement value problem
\begin{eqnarray}\label{div2}
\left\{\begin{array}{l}(\Delta+a^{\alpha}\Delta^{\alpha/2}+b\cdot\nabla+c_k+{\rm div} \hat{b}_k)u_k+f=0\ \ {\rm in}\ D\\
u_k=g\ \ {\rm on}\ D^c
\end{array}\right.\end{eqnarray}
is given by
$$
u_k(x)=E_x\left[e_k(\tau)g(X_{\tau})+\int_0^{\tau}e_k(s)f(X_s)ds\right],\
\ x\in\mathbb{R}^d.
$$

By (\ref{smooth}), (\ref{may2}) and (\ref{may4})--(\ref{13}), we
obtain that
\begin{equation}\label{claim}
\lim\limits_{k\rightarrow\infty}u_k(x)=u^*(x)\ \ {\rm for\ q.e.}\
x\in D.
\end{equation}
We now show  that $u^*$ is a weak solution  to the complement
value problem (\ref{in1}). By (\ref{div1}) and (\ref{13}), we find
that $\{u_k|_D\}$ is a sequence of uniformly bounded functions on
$D$. Let $y\in D$ and $0<r<R$ such that $\overline{B_R(y)}\subset
D$.
 Let $\gamma:[0,\infty)\mapsto[0,\infty)$ be defined as
$$
\left\{\begin{array}{ll}\gamma(s)=1,&\ \ 0\le s\le r,\\
\frac{R-s}{R-r},&\ \ r<s\le R,\\
0,&\ \ s>R
\end{array}\right.$$
and let $\psi(x)=\gamma(|x-y|)$. Then $\phi=\psi^2 u_k\in
W^{1,2}_0(D)$ with compact support in $D$. For $\varepsilon>0$, we
obtain by (\ref{div2}) that
\begin{eqnarray}\label{june1}
&&\int_D\psi^2|\nabla u_k|^2dx\nonumber\\
&=&-2\int_D\psi u_k\langle \nabla u_k,\nabla\psi\rangle dx+\int_D\langle \nabla u_k,\nabla(u_k\psi^2)\rangle dx\nonumber\\
&=&-2\int_D\psi u_k\langle \nabla u_k,\nabla\psi\rangle dx -\frac{a^{\alpha}{\cal A}(d,-\alpha)}{2}\int_{\mathbb{R}^d}\int_{\mathbb{R}^d}\frac{(u_k(x)-u_k(y))((u_k\psi^2)(x)-(u_k\psi^2)(y))}{|x-y|^{d+\alpha}}dxdy\nonumber\\
& &+\int_{D}\langle b,\nabla u_k\rangle u_k\psi^2
dx+\int_{D}c_ku_k^2\psi^2 dx-\int_{D}\langle
\hat{b}_k,\nabla (u_k^2\psi^2)\rangle dx+\int_{D}fu_k\psi^2dx\nonumber\\
&=&-2\int_D\psi u_k\langle \nabla u_k,\nabla\psi\rangle dx\nonumber\\
 & &-\frac{a^{\alpha}{\cal A}(d,-\alpha)}{2}\left(\int_{\mathbb{R}^d}\int_{\mathbb{R}^d}\frac{((u_k\psi)(x)-(u_k\psi)(y))^2}{|x-y|^{d+\alpha}}dxdy
 -\int_{\mathbb{R}^d}\int_{\mathbb{R}^d}\frac{(\psi(x)-\psi(y))^2}{|x-y|^{d+\alpha}}u_k(x)u_k(y)dxdy\right)\nonumber\\
& &+\int_{D}\langle b,\nabla u_k\rangle u_k\psi^2
dx+\int_{D}c_ku_k^2\psi^2 dx-\int_{D}\langle
\hat{b}_k,\nabla (u_k^2\psi^2)\rangle dx+\int_{D}fu_k\psi^2dx\nonumber\\
&\le&\varepsilon\int_D\psi^2|\nabla
u_k|^2dx+\frac{1}{\varepsilon}\int_Du_k^2|\nabla
\psi|^2dx+\frac{a^{\alpha}{\cal A}(d,-\alpha)}{2}\|u_k\|_{\infty}^2\int_{\mathbb{R}^d}\int_{\mathbb{R}^d}\frac{(\psi(x)-\psi(y))^2}{|x-y|^{d+\alpha}}dxdy\nonumber\\
 & &+\frac{\varepsilon}{2}\int_D\psi^2|\nabla
u_k|^2dx+\frac{1}{2\varepsilon}\int_D|b|^2u_k^2\psi^2dx+\int_{D}|c_k|u_k^2\psi^2
dx\nonumber\\
& &+\varepsilon\int_Du_k^2|\nabla
\psi|^2dx+\frac{1}{\varepsilon}\int_D|\hat{b}_k|^2u_k^2\psi^2dx+\varepsilon\int_D\psi^2|\nabla
u_k|^2dx+\frac{1}{\varepsilon}\int_D|\hat{b}_k|^2u_k^2\psi^2dx+\int_{D}|f||u_k|\psi^2dx\nonumber\\
&\le&\frac{5\varepsilon}{2}\int_D\psi^2|\nabla
u_k|^2dx+\frac{\|u_k\|_{\infty}^2}{2\varepsilon}\int_D(2|\nabla
\psi|^2+|b|^2+4|\hat{b}_k|^2)dx\nonumber\\
&&+\|u_k\|_{\infty}^2\left(\int_{D}|c_k|\psi^2
dx+\varepsilon\int_D|\nabla
\psi|^2dx\right)+\|u_k\|_{\infty}\int_{D}|f|dx\nonumber\\
&&+\frac{a^{\alpha}{\cal A}(d,-\alpha)}{2}\|u_k\|_{\infty}^2\int_{\mathbb{R}^d}\int_{\mathbb{R}^d}\frac{(\psi(x)-\psi(y))^2}{|x-y|^{d+\alpha}}dxdy.
\end{eqnarray}
Let $\varepsilon=1/5$. Then, we obtain by the uniform boundedness
of $\{u_k|_D\}$ and (\ref{june1}) that
$$
\sup_{k\in \mathbb{N}}\int_{B_r(y)}|\nabla u_k|^2dx<\infty.
$$
Since $y$ and $r$ are arbitrary, we obtain by (\ref{claim}) that
$u^*|_D\in W^{1,2}_{loc}(D)$.

By taking a subsequence if necessary, we may assume that for any
$\psi\in C_c^\infty(D)$, $u_k\psi\rightarrow u^*\psi$ weakly in
$W^{1,2}(D)$ as $k\rightarrow\infty$ and that
\begin{equation}\label{quasi}u'_k\psi:=\frac{1}{k}\sum_{l=1}^ku_l\psi\rightarrow u^*\psi\
{\rm in}\ W^{1,2}(D)\ {\rm as}\ k\rightarrow\infty.\end{equation}
Let $\phi\in C_c^\infty(D)$. Suppose that ${\rm supp}[\phi]\subset
U\subset {\overline{U}}\subset D$ for some open set $U$. We choose
$\varrho\in C^{\infty}_c(D)$ satisfying $\varrho\equiv1$ on $U$.
Then, there exists a positive constant $C$ which is independent of
$k$ such that
\begin{eqnarray*}
& &\left|\int_{\mathbb{R}^d}\int_{\mathbb{R}^d}\frac{(u'_k(x)-u'_k(y))(\phi(x)-\phi(y))}{|x-y|^{d+\alpha}}dxdy
-\int_{\mathbb{R}^d}\int_{\mathbb{R}^d}\frac{(u(x)-u(y))(\phi(x)-\phi(y))}{|x-y|^{d+\alpha}}dxdy\right|\nonumber\\
&\le&\int_{U}\int_{U}\frac{|((u'_k\varrho-u^*\varrho)(x)-(u'_k\varrho-u^*\varrho)(y))(\phi(x)-\phi(y))|}{|x-y|^{d+\alpha}}dxdy\nonumber\\
& &+\int_{{\rm
supp}[\phi]}\int_{U^c}\frac{|((u'_k-u^*)(x)-(u'_k-u^*)(y))\phi(y)|}{|x-y|^{d+\alpha}}dxdy\nonumber\\
& &+\int_{U^c}\int_{{\rm
supp}[\phi]}\frac{|((u'_k-u^*)(x)-(u'_k-u^*)(y))\phi(x)|}{|x-y|^{d+\alpha}}dxdy\nonumber\\
&\le&C\|u'_k\varrho-u\varrho\|_{W^{1,2}(D)}\|\phi\|_{W^{1,2}(D)}\nonumber\\
&&+\int_{{\rm
supp}[\phi]}\int_{U^c}\frac{2\left(\sup\limits_{k\in\mathbb{N}}\|u'_k\|_{\infty}+\|u^*\|_{\infty}\right)\|\phi\|_{\infty}}{|x-y|^{d+\alpha}}dxdy\nonumber\\
&&+\int_{U^c}\int_{{\rm
supp}[\phi]}\frac{2\left(\sup\limits_{k\in\mathbb{N}}\|u'_k\|_{\infty}+\|u^*\|_{\infty}\right)\|\phi\|_{\infty}}{|x-y|^{d+\alpha}}dxdy.
\end{eqnarray*}
Therefore, we obtain by (\ref{claim}) and the dominated
convergence theorem that
\begin{equation}\label{dom}
\lim_{k\rightarrow\infty}\int_{\mathbb{R}^d}\int_{\mathbb{R}^d}\frac{(u'_k(x)-u'_k(y))(\phi(x)-\phi(y))}{|x-y|^{d+\alpha}}dxdy
=\int_{\mathbb{R}^d}\int_{\mathbb{R}^d}\frac{(u^*(x)-u^*(y))(\phi(x)-\phi(y))}{|x-y|^{d+\alpha}}dxdy.
\end{equation}

Note that
\begin{equation}\label{dom1}
\lim_{k\rightarrow\infty}\int_D\langle\nabla
u'_k,\nabla\phi\rangle dx=\int_D\langle\nabla
u^*,\nabla\phi\rangle dx,
\end{equation}
\begin{equation}\label{dom2}\lim_{k\rightarrow\infty}\int_{D}\langle b,\nabla u'_k\rangle \phi
dx=\int_{D}\langle b,\nabla u^*\rangle \phi
dx,\end{equation}\begin{equation}\label{dom3}
\lim_{k\rightarrow\infty}\int_{D}c_ku'_k\phi dx=\int_{D}cu^*\phi
dx,\end{equation} and
\begin{equation}\label{dom3}\lim_{k\rightarrow\infty}\int_{D}\langle \hat{b}_k,\nabla (u'_k\phi)\rangle
dx=\int_{D}\langle \hat{b},\nabla (u^*\phi)\rangle dx.
\end{equation}
Therefore, we obtain by (\ref{div2}) and (\ref{dom})--(\ref{dom3})
that $u^*$ is a weak solution  to the complement value problem
(\ref{in1}).

\subsection{Continuity of weak solution}\label{dfgh}

In this subsection, we consider the continuity of the weak
solution $u^*$ given by (\ref{weak2}). By Lemma \ref{Holder}, we
know that $u^*|_D$ has a continuous version, which is denoted by
$u|_D$. By (\ref{claim}) and (\ref{quasi}), we find that $u^*|_D$
is quasi-continuous and hence $u(x)=u^*(x)$ for q.e. $x\in D$. We
will show below that if $g$ is continuous at $z\in
\partial D$ then
$$
\lim_{x\rightarrow z}u(x)=g(z).
$$

Note that for $x\in \mathbb{R}^d$,
\begin{eqnarray*}
u^*(x)&=&E_x\left[g(X_{\tau})+\int_0^{\tau}f(X_s)ds\right]\\
&&+E_x\left[(e(\tau)-1)g(X_{\tau})+\int_0^{\tau}(e(s)-1)f(X_s)ds\right].
\end{eqnarray*}
Suppose that $g$ is continuous at $z\in\partial D$. By
\cite[Theorem 1.1]{S}, to prove $\lim_{x\rightarrow z}u(x)=g(z)$,
it suffices to show that there exists an ${\cal E}^0$-exceptional
set $F\subset D$ such that
\begin{equation}\label{Dec1}
\lim_{\stackrel{x\rightarrow z}{x\in D\backslash
F}}E_x[(e(\tau)-1)g(X_{\tau})]=0
\end{equation}
and
\begin{equation}\label{Dec22}
\lim_{\stackrel{x\rightarrow z}{x\in D\backslash
F}}E_x\left[\int_0^{\tau}(e(s)-1)f(X_s)ds\right]=0.
\end{equation}

  For $t>0$, we have
$$E_x[(e(\tau)-1)g(X_{\tau})]
=E_x[(e(\tau)-1)g(X_{\tau});\tau\leq
t]+E_x[(e(\tau)-1)g(X_{\tau});\tau> t].$$ By (\ref{div1}), there
exists an ${\cal E}^0$-exceptional set $F\subset D$ such that for
$x\in D\backslash F$,
\begin{eqnarray}\label{great}
& &\left|E_x[(e(\tau)-1)g(X_{\tau});\tau> t]\right|\nonumber\\
&\leq&\|g\|_\infty\left\{{P}_x(\tau>t)+ E_x[e(\tau);\tau>t]\right\}\nonumber\\
&\le&
\|g\|_\infty\{{P}_x(\tau>t)+(E_x[e^2(\tau)])^{1/2}({P}_x(\tau>t))^{1/2}\}\nonumber\\
&\le&
\|g\|_\infty\{{P}_x(\tau>t)+\Upsilon^{1/8}({P}_x(\tau>t))^{1/2}\}.
\end{eqnarray}
By \cite[Lemma 2.3(3) and Lemma 2.7]{S}, we have
\begin{eqnarray}\label{great21}\limsup_{\stackrel{x\rightarrow z}{x\in
D}}{P}_x(\tau>t)=0.\end{eqnarray} Then, we obtain by (\ref{great})
and (\ref{great21}) that for every $t>0$,
$$\limsup_{\stackrel{x\rightarrow z}{x\in
D\backslash F}}\left|E_x[(e(\tau)-1)g(X_{\tau});\tau>
t]\right|=0.$$ Hence, to prove (\ref{Dec1}), it suffices to show
that
\begin{equation}\label{xcv}
\lim_{t\downarrow0}\sup_{x\in D\backslash
F}|E_x[(e(\tau)-1)g(X_{\tau});\tau\leq t]|=0.
\end{equation}

For $x\in D\backslash F$ and $t > 0$, we obtain by Fatou's lemma
that
\begin{eqnarray*}
&&|E_x[(e(\tau)-1)g(X_{\tau});\tau\leq
t]|\\
&\leq&\|g\|_\infty\liminf_{k\rightarrow\infty}E_x
\left[\left|e^{\int_0^{\tau}(c_k+{\rm div}\,\hat b_k)(X_s)ds}-1\right|;\tau\leq t\right]\\
&\leq&\|g\|_\infty\left\{\sup_{x\in D\backslash F}\limsup_{k\rightarrow\infty}E_x\left[e^{\int_0^{\tau}h_k(X_s)ds}-1;\tau\leq t\right]\right.\\
& &
\left.+\limsup_{k\rightarrow\infty}E_x\left[\left(1-e^{\int_0^{\tau}(c_k+{\rm
div}\,\hat b_k-h_k)(X_s)ds}\right);\tau\leq t\right]\right\}\\
&\leq&\|g\|_\infty\left\{\limsup_{k\rightarrow\infty}E_x\left[e^{\int_0^{t\wedge\tau}h_k(X_s)ds}-1
\right]\right.\\
& & \left.+E_x\left[\left(1-e^{\int_0^{t\wedge\tau}(c_k+{\rm
div}\,\hat b_k-h_k)(X_s)ds}\right)\right]\right\}.
\end{eqnarray*}
By Lemma \ref{lemma23} and Khasminskii's inequality, we get
$$\lim_{t\downarrow0}\sup_{x\in D}\sup_{k\in\mathbb{N}}E_x\left[e^{\int_0^{t\wedge\tau}
h_k(X_s)ds}\right]= 1.$$ Hence, to prove (\ref{xcv}), we need only
show that
$$
\lim_{t\downarrow0}\inf_{x\in
D}\inf_{k\in\mathbb{N}}E_x\left[e^{\int_0^{t\wedge\tau}(c_k+{\rm
div}\,\hat b_k-h_k)(X_s)ds}\right]\ge1.
$$
Further, by Jensen's inequality, we need only show that
$$
\lim_{t\downarrow0}\sup_{x\in
D}\sup_{k\in\mathbb{N}}E_x\left[{\int_0^{t\wedge\tau}(h_k-c_k-{\rm
div}\,\hat b_k)(X_s)ds}\right]=0.
$$
By Lemma \ref{lemma23} and \cite[Theorems 2.1, 2.2 and 4.3]{CH},
we get
\begin{eqnarray*}
& &\sup_{x\in
D}\sup_{k\in\mathbb{N}}E_x\left[{\int_0^{t\wedge\tau}(h_k-c_k-{\rm
div}\,\hat b_k)(X_s)ds}\right]\\
&=&\sup_{x\in D}\sup_{k\in\mathbb{N}}\int_0^t\int_{y\in
D}p^D(s,x,y)(h_k-c_k-{\rm div}\,\hat b_k)(y)dyds\\
&\leq&\sup_{x\in D}\sup_{k\in\mathbb{N}}\int_0^t\int_{y\in
\mathbb{R}^d}p(s,x,y)(h_k-c_k-{\rm div}\,\hat b_k)(y)dyds\\
&\leq&\sup_{x\in D}\sup_{k\in\mathbb{N}}C_1\int_0^t\int_{y\in
\mathbb{R}^d}\left(s^{-d/2}\exp\left(-\frac{C_2|x-y|^2}{s}\right)+s^{-d/2}\wedge\frac{s}{|x-y|^{d+\alpha}}\right)(|h_k|+|c_k|)(y)dyds\\
&&+\sup_{x\in D}\sup_{k\in\mathbb{N}}C_3\int_0^t\int_{y\in
\mathbb{R}^d}\left(s^{-(d+1)/2}\exp\left(-\frac{C_4|x-y|^2}{s}\right)+s^{-(d+1)/2}\wedge\frac{s}{|x-y|^{d+1+\alpha}}\right)|\hat b_k|(y)dyds\\
 &\rightarrow&0\ {\rm as}\ t\downarrow0,
\end{eqnarray*}
where $p^D(t,x,y)$ is the transition density function of the part
process $((X^D_t)_{t\ge 0},(P_x)_{x\in D})$ and $C_i$,
$i=1,2,3,4$, are positive constants. Then (\ref{xcv}) holds and
therefore (\ref{Dec1}) holds.

 For $t>0$, we have
\begin{eqnarray*}E_x\left[\int_0^{\tau}(e(s)-1)f(X_s)ds\right]
&=&E_x\left[\int_0^{\tau}(e(s)-1)f(X_s)ds;\tau\le
t\right]\\
&&+E_x\left[\int_0^{\tau}(e(s)-1)f(X_s)ds;\tau> t\right].
\end{eqnarray*}

By (\ref{12}) and (\ref{13}), we obtain that for $x\in D\backslash
F$,
\begin{eqnarray*}\label{Dec2}
&&\left|E_x\left[\int_0^{\tau}(e(s)-1)f(X_s)ds;\tau> t\right]\right|\\
&\le&\left|E_x\left[\int_0^{\tau}e(s)f(X_s)ds;\tau> t\right]\right|+\left|E_x\left[\int_0^{\tau}f(X_s)ds;\tau> t\right]\right|\\
&\le&\left(E_x\left[\left(\int_0^{\tau}e(s)f(X_s)ds\right)^2\right]\right)^{1/2}(P(\tau>t))^{1/2}\\
&&+(E[\tau;\tau> t])^{1/2}\left(E_x\left[\int_0^{\tau}f^2(X_s)ds\right]\right)^{1/2}\\
&\le&C^{1/2}\Upsilon^{1/8}\left(E_x\left[\tau^4(\tau\vee
1)^4\right]\right)^{1/8}\|f^4\|^{1/4}_{L^{p\vee
1}}(P(\tau>t))^{1/2}\\
&&+C^{1/2}(E[\tau^2])^{1/4}(P(\tau>t))^{1/4}\|f^2\|^{1/2}_{L^{p\vee 1}}\\
&\le&C^{1/2}(\Upsilon^{1/8}+1)\left(E_x\left[\tau^4(\tau\vee
1)^4\right]\right)^{1/8}\|f^4\|^{1/4}_{L^{p\vee
1}}(P(\tau>t))^{1/4}.
\end{eqnarray*}
Then, we obtain by (\ref{15}) and (\ref{great21}) that
$$
\limsup_{\stackrel{x\rightarrow z}{x\in D\backslash
F}}\left|E_x\left[\int_0^{\tau}(e(s)-1)f(X_s)ds;\tau>
t\right]\right|=0.
$$
Hence, to prove (\ref{Dec22}), it suffices to show that
\begin{equation}\label{laole}
\lim_{t\downarrow0}\sup_{x\in D\backslash
F}\left|E_x\left[\int_0^{\tau}(e(s)-1)f(X_s)ds;\tau\le
t\right]\right|=0.
\end{equation}

For $x\in D\backslash F$ and $t > 0$, we obtain by (\ref{12}) and
(\ref{div1}) that
\begin{eqnarray*}
&&\left|E_x\left[\int_0^{\tau}(e(s)-1)f(X_s)ds;\tau\le
t\right]\right|\\
&\leq&E_x\left[\int_0^{t\wedge\tau}e(s)|f|(X_s)ds\right]+E_x\left[\int_0^{t\wedge\tau}|f|(X_s)ds\right]\\
&\le&\left(E_x\left[\int_0^{t\wedge\tau}e^2(s)ds\right]\right)^{1/2}\left(E_x\left[\int_0^{t\wedge\tau}f^2(X_s)ds\right]\right)^{1/2}+E_x\left[\int_0^{t\wedge\tau}|f|(X_s)ds\right]\\
&\le&t^{1/2}\Upsilon^{1/8}\left(E_x\left[\int_0^{\tau}f^2(X_s)ds\right]\right)^{1/2}+t^{1/2}\left(E_x\left[\int_0^{\tau}f^2(X_s)ds\right]\right)^{1/2}\\
&\le&C^{1/2}(\Upsilon^{1/8}+1)\|f^2\|^{1/2}_{L^{p\vee 1}}t^{1/2}.
\end{eqnarray*}
Then (\ref{laole}) holds and therefore (\ref{Dec22}) holds.

\subsection{Uniqueness of solution}

In this subsection, we prove the uniqueness of solution. To this
end, we will show that there exists $M>0$ such that if
$\|h\|_{L^{p\vee 1}}\le M$, then $v\equiv0$ is the unique function
in $B_b(\mathbb{R}^d)$ satisfying $v|_D\in W^{1,2}_{loc}(D)\cap
C(D)$ and
\begin{equation}\label{bl1}
\left\{\begin{array}{l}(\Delta+a^{\alpha}\Delta^{\alpha/2}+b\cdot\nabla+c+{\rm
div}
\hat{b})v=0\ \ {\rm in}\ D,\\
v=0\ \ {\rm on}\ D^c.
\end{array}\right.\end{equation}

Denote by $((X_t)_{t\ge 0}, (P^{\gamma}_x)_{x\in \mathbb{R}^d})$
the Hunt process associated with $({\cal
E}^0_{\gamma},W^{1,2}(\mathbb{R}^d))$. For $\phi\in
W^{1,2}(\mathbb{R}^d)$, we obtain by \cite[Theorem 1.4]{MSW} that
$\phi$ admits a unique Fukushima type decomposition w.r.t.
$((X_t)_{t\ge 0}, (P^{\gamma}_x)_{x\in \mathbb{R}^d})$:
$$ \tilde\phi(X_t)-\tilde\phi(X_0)=M^{\gamma,\phi}_t+N^{\gamma,\phi}_t,\ \
t\ge0,
$$
where $(M^{\gamma,\phi}_t)_{t\ge 0}$ is a locally square
integrable martingale additive functional on $I(\zeta)$ and
$(N^{\gamma,\phi}_t)_{t\ge 0}$ is a continuous additive functional
locally of zero quadratic variation. Hereafter $\zeta$ denotes the
lifetime of $((X_t)_{t\ge 0}, (P^{\gamma}_x)_{x\in \mathbb{R}^d})$
with the totally inaccessible part $\zeta_i$ and
$I(\zeta)=[[0,\zeta[[\cup[[\zeta_i]]$, which is a predictable set
of interval type. Similar to \cite[Lemma 2.2]{CMS1}, we can show
that for $\phi\in C^{\infty}_c(\mathbb{R}^d)$,
\begin{equation}\label{starr}
\int_0^t{\rm div}\phi(X_s)ds=N^{\gamma, \phi^{H}}_t,\ \ t\ge 0.
\end{equation}

Suppose that $v\in B_b(\mathbb{R}^d)$ satisfying $v|_D\in
W^{1,2}_{loc}(D)\cap C(D)$ and (\ref{bl1}). Let
$\{D_n\}_{n\in\mathbb{N}}$ be a sequence of relatively compact
open subsets of $D$ such that $\overline{D}_n\subset D_{n+1}$ and
$D=\cup_{n=1}^{\infty}D_n$, and $\{\chi_n\}_{n\in\mathbb{N}}$  be
a sequence of functions in $C^{\infty}_c(D)$ such that $0\le
\chi_n\le 1$ and $\chi_n|_{D_n}=1$. We have $v\chi_n\in
W^{1,2}_0(D)$. Note that

\begin{equation}\label{san}
\int_{\mathbb{R}^d}\int_{\mathbb{R}^d}\frac{|(v(y)-v(x))(\chi_n(y)-\chi_n(x))|}{|x-y|^{d+\alpha}}dydx<\infty.
\end{equation}
Let $\phi\in W^{1,2}_0(D)$. By (\ref{bl1}) and (\ref{san}), we get
\begin{eqnarray}\label{FOT}
& &{\cal
E}_{\gamma}^0(v\chi_n,\phi)\nonumber\\&=&\int_{\mathbb{R}^d}\langle
\nabla(v\chi_n),\nabla\phi\rangle
dx-\int_{\mathbb{R}^d}\langle b,\nabla (v\chi_n)\rangle \phi dx\nonumber\\
&
&+\frac{a^{\alpha}{\cal A}(d,-\alpha)}{2}\int_{\mathbb{R}^d}\int_{\mathbb{R}^d}\frac{((v\chi_n)(x)-(v\chi_n)(y))(\phi(x)-\phi(y))}{|x-y|^{d+\alpha}}dxdy+(\gamma,v\chi_n\phi)\nonumber\\
&=&{\cal E}^0(v,\chi_n\phi)-\int_{\mathbb{R}^d}(L\chi_n)v\phi
dx-2\int_{\mathbb{R}^d}\langle \nabla v,\nabla\chi_n\rangle\phi
dx\nonumber\\
&
&-a^{\alpha}{\cal A}(d,-\alpha)\int_{\mathbb{R}^d}\left[\int_{\mathbb{R}^d}\frac{(v(y)-v(x))(\chi_n(y)-\chi_n(x))}{|x-y|^{d+\alpha}}dy\right]\phi(x)dx+(\gamma,v\chi_n\phi)\nonumber\\
&=&\left((c+\gamma)v\chi_n-(L\chi_n)v-2\langle \nabla
v,\nabla\chi_n\rangle-a^{\alpha}{\cal A}(d,-\alpha)\int_{\mathbb{R}^d}\frac{(v(y)-v(\cdot))(\chi_n(y)-\chi_n(\cdot))}{|\cdot-y|^{d+\alpha}}dy,\phi\right)\nonumber\\
&&-\int_{\mathbb{R}^d}\langle\hat{b},\nabla(v\chi_n\phi)\rangle dx\nonumber\\
 &:=&(\theta_n,\phi)-\int_{\mathbb{R}^d}\langle\hat{b},\nabla(v\chi_n\phi)\rangle dx.
\end{eqnarray}

Let $n>m$ and $\phi\in C^{\infty}_c(D_m)$. By (\ref{FOT}), we get
\begin{eqnarray}\label{addsept}
& &(\theta_n,\phi)-\int_{\mathbb{R}^d}\langle\hat{b},\nabla(v\chi_n\phi)\rangle dx\nonumber\\
&=&{\cal E}^0(v,\phi)+a^{\alpha}{\cal A}(d,-\alpha)\int_{D_m}\int_{D^c_n}\frac{v(y)(1-\chi_n(y))}{|x-y|^{d+\alpha}}dy\phi(x) dx+(\gamma,v\phi)\nonumber\\
&=&\left((c+\gamma)v+a^{\alpha}{\cal A}(d,-\alpha)\int_{D\cap
D^c_n}\frac{v(y)(1-\chi_n(y))}{|\cdot-y|^{d+\alpha}}dy,\phi\right)-\int_{\mathbb{R}^d}\langle\hat{b},\nabla(v\phi)\rangle dx.
\end{eqnarray}
Since $\phi\in C^{\infty}_c(D_m)$ is arbitrary, by (\ref{addsept}), we find that for $n>m$,
$$
\theta_n(x)=(c(x)+\gamma)v(x)+a^{\alpha}{\cal
A}(d,-\alpha)\int_{D\cap
D^c_n}\frac{v(y)(1-\chi_n(y))}{|x-y|^{d+\alpha}}dy,\ \ x\in D_m.
$$
Then,
\begin{equation}\label{nm2}\theta_n\ {\rm converges\ to}\ (c+\gamma)v\ {\rm uniformly\ on\ any\ compact\ subset\
of}\ D\ {\rm as}\ n\rightarrow\infty.
\end{equation}

We fix a function $\xi\in C(\mathbb{R}^d)\cap
L^1(\mathbb{R}^d;dx)$ with $0<\xi\le 1$. Set ${\cal
G}=G^{\gamma}_1\xi$ and $\hat{\cal G}=\hat{G}^{\gamma}_1\xi$,
where $(G^{\gamma}_{\iota})_{\iota>0}$ and
$(\hat{G}^{\gamma}_{\iota})_{\iota>0}$ denote the resolvent and
co-resolvent associated with $({\cal
E}_{\gamma},W^{1,2}(\mathbb{R}^d))$, respectively. Define
\begin{eqnarray*}
\Theta&:=&\{\{{V_n}\}\,|\, {V_n}\ \mbox{is}\ {\cal
E}_{\gamma}^0{\textrm{-}}\mbox{quasi} {\textrm{-}}\mbox{open},\
{V_n}\subset V_{n+1}\ {\cal E}^0_{\gamma}
{\textrm{-}{\rm q.e.}} \\
         &&\ \ \ \ \ \ \ \ \ \ \ \forall~n\in \mathbb{N},\
         \mbox{and}\ E=\bigcup_{n=1}^{\infty}{V_n}\ {\cal E}_{\gamma}^0{\textrm{-}{\rm
         q.e.}}\}.
\end{eqnarray*}
By \cite{CMS,MSW}, there exists a $\{V_l\}\in \Theta$ such that
$\tilde{\hat{\cal G}}$, $\widetilde{\hat{b}^{H}}$ and
$\widetilde{\hat{b}_k^{H}}$, $k\in \mathbb{N}$, are all bounded on
each $V_l$. By \cite[ptoposition 2.18]{MOR}, we may assume without
loss of generality that there exists a sequence $\{\phi_n\}\subset
C^{\infty}_c(\mathbb{R}^d)$ satisfying
\begin{equation}\label{nmnm}
\phi_n\rightarrow \hat{b}^{H}\ {\rm  in}\ W^{1,2}(\mathbb{R}^d)\
{\rm  and}\ \phi_n\rightarrow\widetilde{\hat{b}^{H}}\ {\rm
uniformly\ on\ each}\ V_l\ {\rm as}\ n\rightarrow\infty.
\end{equation}
Let $V$ be a quasi-open set of $\mathbb{R}^d$. We denote
$W^{1,2}_{0}(V):=\{\phi\in W^{1,2}(\mathbb{R}^d):\phi=0\ {\rm on}\
V^c\}$. For $l\in \mathbb{N}$, we define
$E_l=\{x\in\mathbb{R}^d:\widetilde{G_1^{\gamma,V_l}\xi}(x)>1/l\}$,
where $(\hat{G}^{\gamma,V_l}_{\iota})_{\iota>0}$ is the
co-resolvent associated with the part Dirichlet form $({\cal
E}^0_{\gamma},W^{1,2}_0(V_l))$.  Then $\{E_l\}\in \Theta$
satisfying $\overline{E}^{{\cal E}^0_{\gamma}}_l\subset E_{l+1}$
q.e. and $E_l\subset V_l$ q.e. for each $l\in \mathbb{N}$. Here
$\overline{E}^{{\cal E}^0_{\gamma}}_l$ denotes the ${\cal
E}_{\gamma}^0$-quasi-closure of $E_l$. Define
$\omega_l=l\widetilde{G_1^{\gamma,V_l}\xi}\wedge 1$. Then
$\omega_l\in B_b(\mathbb{R}^d)\cap W^{1,2}(\mathbb{R}^d)$,
$\omega_l=1$ on $E_l$ and $\omega_l=0$ on $V_l^c$. For
$\phi,\psi\in W^{1,2}_{loc}(\mathbb{R}^d)$, we define the
stochastic integral $\int_0^t\phi(X_{s-})dA^{\psi}_s$ as in
\cite{CMS}.

By \cite[Theorems 1.4 and 1.15]{MSW},
$(M^{\gamma,\hat{b}^{H}}_{t\wedge\tau_{E_l}})_{t\ge 0}$ is a
$P^{\gamma}_{x}$-square-integrable martingale for q.e. $x\in
\mathbb{R}^d$ and $l\in \mathbb{N}$. For $\phi\in
W^{1,2}(\mathbb{R}^d)$, we define $$
A^{\phi}_t=\tilde\phi(X_t)-\tilde\phi(X_0), \ \ t\ge 0.
$$
For $l\in\mathbb{N}$, $n>m$, $t\ge0$ and $x\in E_l\cap D_m$, we
define
$$
B^{l,m,n}_t:=\int_0^{t\wedge\tau_{E_l\cap
D_m}}\theta_n(X_s)ds+\int_0^{t\wedge\tau_{E_l\cap
D_m}}v(X_{s-})dA^{\hat{b}^{H}}_s$$ and $$
c^{l,m,n}_t(x):=E^{\gamma}_x[B^{l,m,n}_t].
$$
Denote by $(p^{\gamma,E_l\cap D_m}(t,x,y))_{t\ge 0}$ and
$(p_t^{\gamma,E_l\cap D_m})_{t\ge 0}$ the transition density
function and the transition semigroup of the part process
$((X^{E_l\cap D_m}_t)_{t\ge 0},(P^{\gamma}_x)_{x\in E_l\cap
D_m})$, respectively. By the continuity of the function $t\mapsto
p^{\gamma,E_l\cap D_m}(t,x,y)$, which can be proved similar to
\cite[Lemma 2.3(4)]{S}, we know that the function $t\mapsto
c^{l,m,n}_t(x)$ is continuous for $dx$-a.e. $x\in E_l\cap D_m$. We
have $c^{l,m,n}_t\in L^2(E_l\cap D_m;dx)$ for $t\ge 0$ and
\begin{equation}\label{nm4}
c^{l,m,n}_{t+s}(x)=c^{l,m,n}_t(x)+p_t^{\gamma,E_l\cap
D_m}c^{l,m,n}_s(x),\ \ t,s\ge 0.
\end{equation}

Denote by $(\hat{T}_t^{\gamma,E_l\cap D_m})_{t\ge 0}$ and
$(\hat{G}_{\iota}^{\gamma,E_l\cap D_m})_{\iota>0}$  the
co-semigroup and the co-resolvent of $({\cal
E}_{\gamma}^0,$ $W^{1,2}_{0} (E_l\cap D_m))$ on $L^2(E_l\cap
D_m;dx)$, respectively. Let $\phi\in \hat{G}_1^{\gamma,E_l\cap
D_m}(B_b(E_l\cap D_m))$. Since $\tilde{\hat{\cal G}}$ is  bounded
on $V_l$, both $\phi$ and $\hat{G}_1^{\gamma,E_l\cap D_m}|\phi|$
are bounded on $E_l\cap D_m$. Hence $\hat{T}_s^{\gamma,E_l\cap
D_m}\hat{G}_1^{\gamma,E_l\cap D_m}|\phi|$ is bounded on $E_l\cap
D_m$ and $v\hat{T}_s^{\gamma,E_l\cap D_m}\hat{G}_1^{\gamma,E_l\cap
D_m}\phi\in W^{1,2}(\mathbb{R}^d)$ for all $s\ge0$. Further, we
have
\begin{eqnarray}\label{coer0}
&&\sup_{0\le s\le t}\|\hat{T}^{\gamma,E_l\cap
D_m}_s\hat{G}_1^{\gamma,E_l\cap D_m}\phi\|_{\infty}\nonumber\\
&\le&e^t\|\hat{G}_1^{\gamma,E_l\cap
D_m}|\phi|\|_{\infty}\nonumber\\
&<&\infty,
\end{eqnarray}
\begin{eqnarray}\label{coer1}
&&\lim_{t\downarrow0}\sup_{0\le s\le t}\|\hat{T}^{\gamma,E_l\cap
D_m}_s\hat{G}_1^{\gamma,E_l\cap D_m}\phi-\hat{G}_1^{\gamma,E_l\cap
D_m}\phi\|_{\infty}\nonumber\\
&=&\lim_{t\downarrow0}\sup_{0\le s\le
t}\|e^{-s}\hat{T}^{\gamma,E_l\cap D_m}_s\hat{G}_1^{\gamma,E_l\cap
D_m}\phi-\hat{G}_1^{\gamma,E_l\cap
D_m}\phi\|_{\infty}\nonumber\\
&=&\lim_{t\downarrow0}\sup_{0\le s\le
t}\left\|\int_0^se^{-w}\hat{T}^{\gamma,E_l\cap
D_m}_w\phi dw\right\|_{\infty}\nonumber\\
&\le&\lim_{t\downarrow0}t\|\phi\|_{\infty}\nonumber\\
&=&0,
\end{eqnarray}
\begin{eqnarray}\label{coer02}
&&\sup_{0\le s\le t}{\cal E}^0_{\gamma+1}(\hat{T}^{\gamma,E_l\cap
D_m}_s\hat{G}_1^{\gamma,E_l\cap D_m}\phi,\hat{T}^{\gamma,E_l\cap
D_m}_s\hat{G}_1^{\gamma,E_l\cap D_m}\phi)\nonumber\\
&=&\sup_{0\le s\le t}(\hat{T}^{\gamma,E_l\cap
D_m}_s\hat{G}_1^{\gamma,E_l\cap D_m}\phi,\hat{T}^{\gamma,E_l\cap
D_m}_s\phi)\nonumber\\
&<&\infty,
\end{eqnarray}
and
\begin{eqnarray}\label{coer2}
&&\lim_{t\downarrow0}\sup_{0\le s\le t}{\cal
E}^0_{\gamma+1}(\hat{T}^{\gamma,E_l\cap
D_m}_s\hat{G}_1^{\gamma,E_l\cap D_m}\phi-\hat{G}_1^{\gamma,E_l\cap
D_m}\phi,\hat{T}^{\gamma,E_l\cap D_m}_s\hat{G}_1^{\gamma,E_l\cap
D_m}\phi-\hat{G}_1^{\gamma,E_l\cap
D_m}\phi)\nonumber\\
&=&\lim_{t\downarrow0}\sup_{0\le s\le t}(\hat{T}^{\gamma,E_l\cap
D_m}_s\hat{G}_1^{\gamma,E_l\cap D_m}\phi-\hat{G}_1^{\gamma,E_l\cap
D_m}\phi,\hat{T}^{\gamma,E_l\cap
D_m}_s\phi-\phi)\nonumber\\
&=&0.
\end{eqnarray}

Denote
$L^{\gamma}:=\Delta+a^{\alpha}\Delta^{\alpha/2}+(b-\hat{b})\cdot\nabla-\gamma$.
By \cite[Proposition 4.1.10]{CF} (note that the assertion holds
true also in the setting of semi-Dirichlet forms), \cite[(22) and
Theorem 2.8 with its proof]{CMS}, (\ref{nmnm}), (\ref{coer0}),
(\ref{coer02}), \cite[I. Corollary 4.15]{MR} and the equivalence
of the ${\cal E}^0_{\gamma}$-norm and the $W^{1,2}$-norm, we get
\begin{eqnarray}\label{CF}
&&E^{\gamma}_{\hat{G}_1^{\gamma,E_l\cap D_m}\phi\cdot
dx}\left[\int_0^{t\wedge\tau_{E_l\cap
D_m}}v(X_{s-})dA^{\hat{b}^{H}}_s\right]\nonumber\\
&=&\lim_{n\rightarrow\infty}E^{\gamma}_{\hat{G}_1^{\gamma,E_l\cap
D_m}\phi\cdot dx}\left[\int_0^{t\wedge\tau_{E_l\cap
D_m}}v(X_{s-})dA^{\phi_n}_s\right]\nonumber\\
&=&\lim_{n\rightarrow\infty}E^{\gamma}_{\hat{G}_1^{\gamma,E_l\cap
D_m}\phi\cdot dx}\left[\int_0^{t\wedge\tau_{E_l\cap
D_m}}v(X_{s})L^{\gamma}\phi_n(X_s)ds\right]\nonumber\\
&=&-\lim_{n\rightarrow\infty}\int_0^t{\cal
E}^0_{\gamma}(\phi_n,v\hat{T}^{\gamma,E_l\cap
D_m}_s\hat{G}_1^{\gamma,E_l\cap
D_m}\phi)ds\nonumber\\
&=&-\int_0^t{\cal
E}^0_{\gamma}(\hat{b}^{H},v\hat{T}^{\gamma,E_l\cap
D_m}_s\hat{G}_1^{\gamma,E_l\cap D_m}\phi)ds.
\end{eqnarray}
Then, we obtain by (\ref{bbb}), (\ref{FOT}), (\ref{coer1}),
(\ref{coer2}), (\ref{CF}), \cite[I. Corollary 4.15]{MR} and the
equivalence of the ${\cal E}^0_{\gamma}$-norm and the
$W^{1,2}$-norm that
\begin{eqnarray}\label{nm5}
&&\lim_{t\rightarrow 0}\frac{1}{t}E^{\gamma}_{\hat{G}_1^{\gamma,E_l\cap D_m}\phi\cdot dx}[B^{l,m,n}_t]\nonumber\\
&=&\lim_{t\rightarrow 0}\frac{1}{t}\int_0^t(p_s^{\gamma,{E_l\cap
D_m}}\theta_n,\hat{G}_1^{\gamma,E_l\cap
D_m}\phi)ds+\lim_{t\rightarrow
0}\frac{1}{t}E^{\gamma}_{\hat{G}_1^{\gamma,E_l\cap D_m}\phi\cdot
dx}\left[\int_0^{t\wedge\tau_{E_l\cap
D_m}}v(X_{s-})dA^{\hat{b}^{H}}_s\right]\nonumber\\
&=&(\theta_n,\hat{G}_1^{\gamma,E_l\cap D_m}\phi)-{\cal E}^0_{\gamma}(\hat{b}^{H},v\hat{G}_1^{\gamma,E_l\cap D_m}\phi)\nonumber\\
&=&(\theta_n,\hat{G}_1^{\gamma,E_l\cap D_m}\phi)-\int_{\mathbb{R}^d}\langle\hat{b},\nabla(v\hat{G}_1^{\gamma,E_l\cap D_m}\phi)\rangle dx\nonumber\\
&=&{\cal E}^0_{\gamma}(v\chi_n,\hat{G}_1^{\gamma,E_l\cap
D_m}\phi).
\end{eqnarray}

Define
$$
\eta_{l,m,n}(x)=E^{\gamma}_x[(v\chi_n)(X_{\tau_{E_l\cap D_m}})],\
\ x\in\mathbb{R}^d.
$$
We have
\begin{equation}\label{reduced}
\eta_{l,m,n}(x)=E^{\gamma}_x[\eta_{l,m,n}(X_{t\wedge\tau_{E_l\cap
D_m}})], \ \ t\ge 0,\ x\in E_l\cap D_m,
\end{equation} and $\eta_{l,m,n}(x)=v\chi_n(x)$ for
q.e. $x\in (E_l\cap D_m)^c$. By \cite[Theorem 3.5.1]{Oshima}, we
get
\begin{equation}\label{wanle}
{\cal E}^0_{\gamma}(v\chi_n,\hat{G}_1^{\gamma,E_l\cap
D_m}\phi)={\cal
E}^0_{\gamma}(v\chi_n-\eta_{l,m,n},\hat{G}_1^{\gamma,E_l\cap
D_m}\phi), \ \ \forall \phi\in \hat{G}_1^{\gamma,E_l\cap
D_m}(B_b(E_l\cap D_m)).
\end{equation}
Define
\begin{equation}\label{eas}
\hat{S}^{l,m}_t:=\int_0^t\hat{T}^{\gamma,E_l\cap D_m}_sds,\ \ t\ge
0.
\end{equation}
Similar to \cite[(1.5.5), page 39]{Fuku}, we can show
that\begin{equation}\label{nm6} {\cal
E}_{\gamma}^0(\varpi,\hat{S}^{l,m}_t\rho)=(\varpi,\rho-\hat{T}_t^{\beta,D_m}\rho),\
\ \forall \varpi\in W^{1,2}_{0}(E_l\cap D_m),\ \rho\in L^2(E_l\cap
D_m;dx).
\end{equation}
Then, we obtain by (\ref{nm4}), (\ref{nm5}) and
(\ref{wanle})--(\ref{nm6}) that for $\phi\in
\hat{G}_1^{\gamma,E_l\cap D_m}(B_b(E_l\cap D_m))$ and $t,r>0$,
\begin{eqnarray*}
&&(c^{l,m,n}_t,\hat{G}_1^{\gamma,E_l\cap
D_m}\phi-\hat{T}^{\gamma,E_l\cap D_m}_r\hat{G}_1^{\gamma,E_l\cap
D_m}\phi)\\&=&\lim_{s\rightarrow
0}\frac{1}{s}(c^{l,m,n}_t,\hat{S}^{l,m}_r\hat{G}_1^{\gamma,E_l\cap
D_m}\phi-\hat{T}^{\gamma,E_l\cap
D_m}_s\hat{S}^{l,m}_r\hat{G}_1^{\gamma,E_l\cap
D_m}\phi)\\
&=&\lim_{s\rightarrow
0}\frac{1}{s}(c^{l,m,n}_s,\hat{S}^{l,m}_r\hat{G}_1^{\gamma,E_l\cap
D_m}\phi-\hat{T}^{\gamma,E_l\cap
D_m}_t\hat{S}^{l,m}_r\hat{G}_1^{\gamma,E_l\cap
D_m}\phi)\\
&=&{\cal
E}^0_{\gamma}(v\chi_n,\hat{S}^{l,m}_r\hat{G}_1^{\gamma,E_l\cap
D_m}\phi-\hat{T}^{\gamma,E_l\cap
D_m}_t\hat{S}^{l,m}_r\hat{G}_1^{\gamma,E_l\cap
D_m}\phi)\\
&=&{\cal
E}^0_{\gamma}(v\chi_n-\eta_{l,m,n},\hat{S}^{l,m}_r\hat{G}_1^{\gamma,E_l\cap
D_m}\phi-\hat{T}^{\gamma,E_l\cap
D_m}_t\hat{S}^{l,m}_r\hat{G}_1^{\gamma,E_l\cap
D_m}\phi)\\
&=&(v\chi_n-\eta_{l,m,n},\hat{G}_1^{\gamma,E_l\cap
D_m}\phi-\hat{T}^{\gamma,E_l\cap D_m}_t\hat{G}_1^{\gamma,E_l\cap
D_m}\phi\\
&&\ \ \ \ -\hat{T}^{\gamma,E_l\cap D_m}_r\hat{G}_1^{\gamma,E_l\cap
D_m}\phi+\hat{T}^{\gamma,E_l\cap
D_m}_{t+r}\hat{G}_1^{\gamma,E_l\cap
D_m}\phi)\\
&=&(v\chi_n-\eta_{l,m,n}-p^{\gamma,E_l\cap
D_m}_t(v\chi_n-\eta_{l,m,n}),\hat{G}_1^{\gamma,E_l\cap
D_m}\phi-\hat{T}^{\gamma,E_l\cap D_m}_r\hat{G}_1^{\gamma,E_l\cap
D_m}\phi).
\end{eqnarray*}
Hence
$\kappa^{l,m,n}_t:=(c^{l,m,n}_t-(v\chi_n-\eta_{l,m,n})+p^{\gamma,E_l\cap
D_m}_t(v\chi_n-\eta_{l,m,n}),\hat{G}_1^{\gamma,E_l\cap D_m}\phi)$
satisfies the linear equation
$\kappa^{l,m,n}_t=\kappa^{l,m,n}_{t+r}-\kappa^{l,m,n}_r$. By
(\ref{nm5}) and (\ref{wanle}), we get $\lim_{t\rightarrow
0}\kappa^{l,m,n}_t/t=0$. Then, $\kappa^{l,m,n}_t=0$. Since
$\phi\in \hat{G}_1^{\gamma,E_l\cap D_m}(B_b(E_l\cap D_m))$ is
arbitrary, we obtain by the continuity of the function $t\mapsto
p^{\gamma,E_l\cap D_m}(t,x,y)$ and the continuity of the function
$t\mapsto c^{l,m,n}_t(x)$ that for $dx$-a.e. $x\in E_l\cap D_m$,
\begin{eqnarray*}
(v\chi_n-\eta_{l,m,n})(x)&=&E^{\gamma}_x[(v\chi_n-\eta_{l,m,n})(X_{t\wedge\tau_{E_l\cap
D_m}})]+E^{\gamma}_x\left[\int_0^{t\wedge\tau_{E_l\cap
D_m}}\theta_n(X_s)ds\right]\\
& &+E^{\gamma}_x\left[\int_0^{t\wedge\tau_{E_l\cap
D_m}}v(X_{s-})dA^{\hat{b}^{H}}_s\right]\\
&=&E^{\gamma}_x[(v\chi_n-\eta_{l,m,n})(X_{t\wedge\tau_{E_l\cap
D_m}})]+E^{\gamma}_x\left[\int_0^{t\wedge\tau_{E_l\cap
D_m}}\theta_n(X_s)ds\right]\\
& &+E^{\gamma}_x\left[\int_0^{t\wedge\tau_{E_l\cap
D_m}}v(X_{s-})dN^{\gamma,\hat{b}^{H}}_s\right], \ \ \forall t\ge
0.
\end{eqnarray*} By (\ref{reduced}), we obtain that for $dx$-a.e.
$x\in E_l\cap D_m$,
\begin{eqnarray}\label{zl}
(v\chi_n)(x)&=&E^{\gamma}_x[(v\chi_n)(X_{t\wedge\tau_{E_l\cap
D_m}})]+E^{\gamma}_x\left[\int_0^{t\wedge\tau_{E_l\cap
D_m}}\theta_n(X_s)ds\right]\nonumber\\
& &+E^{\gamma}_x\left[\int_0^{t\wedge\tau_{E_l\cap
D_m}}v(X_{s-})dN^{\gamma,\hat{b}^{H}}_s\right], \ \ \forall t\ge
0.
\end{eqnarray}

 Note that $v\in
B_b(\mathbb{R}^d)$ and $v=0$ on $D^c$. Letting
$n\rightarrow\infty$, we obtain by (\ref{nm2}) and (\ref{zl}) that
for $dx$-a.e. $x\in E_l\cap D_m$,
\begin{eqnarray}\label{r2}
v(x)&=&E^{\gamma}_x[v(X_{t\wedge\tau_{E_l\cap
D_m}})]+E^{\gamma}_x\left[\int_0^{t\wedge\tau_{E_l\cap
D_m}}((c+\gamma)v)(X_s)ds\right]\nonumber\\
& &+E^{\gamma}_x\left[\int_0^{t\wedge\tau_{E_l\cap
D_m}}v(X_{s-})dN^{\gamma,\hat{b}^{H}}_s\right], \ \ \forall t\ge
0.
\end{eqnarray}
Letting $m\rightarrow\infty$, we obtain by (\ref{r2}) that for
$dx$-a.e. $x\in E_l\cap D$,
\begin{eqnarray}\label{chan}
v(x)&=&E^{\gamma}_x[v(X_{t\wedge\tau_{E_l\cap
D}})]+E^{\gamma}_x\left[\int_0^{t\wedge\tau_{E_l\cap
D}}((c+\gamma)v)(X_s)ds\right]\nonumber\\
& &+E^{\gamma}_x\left[\int_0^{t\wedge\tau_{E_l\cap
D}}v(X_{s-})dN^{\gamma,\hat{b}^{H}}_s\right], \ \ \forall t\ge 0.
\end{eqnarray}

Define \begin{equation}\label{ccc} {\cal
I}^l_t=v(X_t)1_{\{\tau_{E_l\cap D}>t\}}+\int_0^{t\wedge
\tau_{E_l\cap D}}((c+\gamma)v)(X_s)ds+\int_0^{t\wedge\tau_{E_l\cap
D}}v(X_{s-})dN^{\gamma,\hat{b}^{H}}_s.
\end{equation}
 By (\ref{chan}), we find that $({\cal I}^l_t)_{t\ge 0}$ is a locally bounded martingale under $P^{\gamma}_x$  for $dx$-a.e. $x\in E_l\cap D$.
For $k\in \mathbb{N}$ and $t\ge 0$, define
$$
e^{k}_{\gamma}(t):=e^{\int_0^t(c_k+{\rm div}\,\hat
b_k+\gamma)(X_s)ds}.
$$
 The integration by parts formula for semi-martingales implies
that
$$
e^{k}_{\gamma}(t){\cal I}^l_t-v(x)=\int_0^t{\cal
I}^l_sde^{k}_{\gamma}(s)+\int_0^te^{k}_{\gamma}(s)d{\cal I}^l_s.
$$
By (\ref{starr}), we have
\begin{equation}\label{swim}
N^{\gamma,\phi_n-\hat{b}_k^{H}}_t=\int_0^t(L^{\gamma}\phi_n-{\rm
div}\,\hat b_k)(X_s)ds,\ \ t\ge 0.
\end{equation}
By (\ref{nmnm}) (\ref{ccc}), (\ref{swim}), \cite[Theorem 2.8 and
(22)]{CMS} and \cite[(11) and (12)]{MSW}, we obtain that for
$dx$-a.e. $x\in E_l\cap D$ and $t\ge 0$, \begin{eqnarray}\label{Oshi} &&v(x)+\int_0^te^{k}_{\gamma}(s)d{\cal I}^l_s\nonumber\\
&=& e^{k}_{\gamma}(t){\cal I}^l_t-\int_0^t{\cal
I}^l_sde^{k}_{\gamma}(s)\nonumber\\
&=&e^{k}_{\gamma}(t)v(X_t)1_{\{\tau_{E_l\cap
D}>t\}}+e^{k}_{\gamma}(t)\int_0^{t\wedge
\tau_{E_l\cap D}}((c+\gamma)v)(X_s)ds\nonumber\\
&&+e^{k}_{\gamma}(t)\left(\lim_{n\rightarrow\infty}\int_0^{t\wedge\tau_{E_l\cap D}}v(X_{s-})dN^{\gamma,\phi_n}_s\right)\nonumber\\
&&-\int_0^tv(X_s)1_{\{\tau_{E_l\cap
D}>s\}}de^{k}_{\gamma}(s)-\int_0^t\int_0^{s\wedge
\tau_{E_l\cap D}}((c+\gamma)v)(X_w)dwde^{k}_{\gamma}(s)\nonumber\\
& &-\lim_{n\rightarrow\infty}\int_0^t\int_0^{s\wedge\tau_{E_l\cap
D}}v(X_{w-})dN^{\gamma,\phi_n}_wde^{k}_{\gamma}(s)\nonumber\\
&=&e^{k}_{\gamma}(t)v(X_t)1_{\{\tau_{E_l\cap
D}>t\}}+e^{k}_{\gamma}(t)\int_0^{t\wedge
\tau_{E_l\cap D}}((c+\gamma)v)(X_s)ds\nonumber\\
&&+e^{k}_{\gamma}(t)\left(\lim_{n\rightarrow\infty}\int_0^{t\wedge\tau_{E_l\cap D}}v(X_{s-})dN^{\gamma,\phi_n}_s\right)\nonumber\\
&&-\int_0^te^{k}_{\gamma}(s)((c_k+{\rm div}\,\hat
b_k+\gamma)v)(X_s)1_{\{\tau_{E_l\cap D}>
s\}}ds\nonumber\\
& &-\int_0^t((c+\gamma)v)(X_w)1_{\{\tau_{E_l\cap D}\ge
w\}}\left(\int_w^t e^{k}_{\gamma}(s)(c_k+{\rm div}\,\hat
b_k+\gamma)(X_s)ds\right)dw\nonumber\\
& &-\lim_{n\rightarrow\infty}\int_0^tv(X_{w-})1_{\{\tau_{E_l\cap
D}\ge w\}}\left(\int_w^t e^{k}_{\gamma}(s)(c_k+{\rm div}\,\hat
b_k+\gamma)(X_s)ds\right)dN^{\gamma,\phi_n}_w\nonumber\\
&=&e^{k}_{\gamma}(t)v(X_t)1_{\{\tau_{E_l\cap
D}>t\}}+\int_0^te^{k}_{\gamma}(s)((c-c_k-{\rm div}\,\hat
b_k)v)(X_s)1_{\{\tau_{E_l\cap D}>
s\}}ds\nonumber\\
&&+\lim_{n\rightarrow\infty}\int_0^te^{k}_{\gamma}(s)v(X_{s-})1_{\{\tau_{E_l\cap
D}>
s\}}dN^{\gamma,\phi_n}_s\nonumber\\
&=&e^{k}_{\gamma}(t)v(X_t)1_{\{\tau_{E_l\cap
D}>t\}}+\int_0^te^{k}_{\gamma}(s)((c-c_k)v)(X_s)1_{\{\tau_{E_l\cap
D}>
s\}}ds\nonumber\\
&&+\lim_{n\rightarrow\infty}\int_0^te^{k}_{\gamma}(s)v(X_{s-})1_{\{\tau_{E_l\cap
D}>
s\}}dN^{\gamma,\phi_n-\hat{b}_k^{H}}_s\nonumber\\
&=&e^{k}_{\gamma}(t)v(X_t)1_{\{\tau_{E_l\cap
D}>t\}}+\int_0^{t\wedge\tau_{E_l\cap D}}e^{k}_{\gamma}(s)((c-c_k)v)(X_s)ds\nonumber\\
&&+\lim_{n\rightarrow\infty}e^{k}_{\gamma}(t\wedge \tau_{E_l\cap D})\int_0^{t\wedge\tau_{E_l\cap D}}v(X_{s-})dN^{\gamma,\phi_n-\hat{b}_k^{H}}_s\nonumber\\
&&-\lim_{n\rightarrow\infty}\int_0^{t\wedge\tau_{E_l\cap
D}}\left(\int_0^sv(X_{w-})dN^{\gamma,\phi_n-\hat{b}_k^{H}}_w\right)e^{k}_{\gamma}(s)(c_k+{\rm
div}\,\hat
b_k+\gamma)(X_s)ds\nonumber\\
&=&e^{k}_{\gamma}(t)v(X_t)1_{\{\tau_{E_l\cap
D}>t\}}+\int_0^{t\wedge\tau_{E_l\cap D}}e^{k}_{\gamma}(s)((c-c_k)v)(X_s)ds\nonumber\\
&&+\lim_{n\rightarrow\infty}e^{k}_{\gamma}(t\wedge \tau_{E_l\cap D})\int_0^{t\wedge\tau_{E_l\cap D}}v(X_{s-})dN^{\gamma,\phi_n-\hat{b}_k^{H}}_s\nonumber\\
&&+\lim_{n\rightarrow\infty}\int_0^{t\wedge\tau_{E_l\cap
D}}\left(\int_0^sv(X_{w-})dN^{\gamma,\phi_n-\hat{b}_k^{H}}_w\right)e^{k}_{\gamma}(s)(h_k-c_k-{\rm
div}\,\hat
b_k)(X_s)ds\nonumber\\
&&-\lim_{n\rightarrow\infty}\int_0^{t\wedge\tau_{E_l\cap
D}}\left(\int_0^sv(X_{w-})dN^{\gamma,\phi_n-\hat{b}_k^{H}}_w\right)e^{k}_{\gamma}(s)(h_k+\gamma)(X_s)ds \nonumber\\
&=&e^{k}_{\gamma}(t)v(X_t)1_{\{\tau_{E_l\cap
D}>t\}}+\int_0^{t\wedge\tau_{E_l\cap D}}e^{k}_{\gamma}(s)((c-c_k)v)(X_s)ds\nonumber\\
&&+e^{k}_{\gamma}(t\wedge \tau_{E_l\cap D})\int_0^{t\wedge\tau_{E_l\cap D}}v(X_{s-})dN^{\gamma,\hat{b}^{H}-\hat{b}_k^{H}}_s\nonumber\\
&&+\int_0^{t\wedge\tau_{E_l\cap
D}}\left(\int_0^sv(X_{w-})dN^{\gamma,\hat{b}^{H}-\hat{b}_k^{H}}_w\right)e^{k}_{\gamma}(s)(h_k-c_k-{\rm
div}\,\hat
b_k)(X_s)ds\nonumber\\
&&-\int_0^{t\wedge\tau_{E_l\cap
D}}\left(\int_0^sv(X_{w-})dN^{\gamma,\hat{b}^{H}-\hat{b}_k^{H}}_w\right)e^{k}_{\gamma}(s)(h_k+\gamma)(X_s)ds
.
\end{eqnarray}

Note that
\begin{eqnarray}\label{shui}
&&\left|\int_0^{t\wedge\tau_{E_l\cap
D}}\left(\int_0^sv(X_{w-})dN^{\gamma,\hat{b}^{H}-\hat{b}_k^{H}}_w\right)e^{k}_{\gamma}(s)(h_k-c_k-{\rm
div}\,\hat
b_k)(X_s)ds\right|\nonumber\\
&&+\left|\int_0^{t\wedge\tau_{E_l\cap
D}}\left(\int_0^sv(X_{w-})dN^{\gamma,\hat{b}^{H}-\hat{b}_k^{H}}_w\right)e^{k}_{\gamma}(s)(h_k+\gamma)(X_s)ds\right|\nonumber\\
&\le&\left(\sup_{0\le s\le t\wedge\tau_{E_l\cap D}}\left|\int_0^sv(X_{w-})dN^{\gamma,\hat{b}^{H}-\hat{b}_k^{H}}_w\right|\right)e^{\int_0^{t\wedge\tau_{E_l\cap D}}(h_k+\gamma)(X_s)ds}\nonumber\\
&&\ \ \ \ \cdot\left(\int_0^{t\wedge\tau_{E_l\cap
D}}(2h_k-c_k+\gamma)(X_s)ds-N^{\gamma,\hat{b}_k^{H}}_{t\wedge\tau_{E_l\cap
D}}\right).
\end{eqnarray}
Define
$$
e_{\gamma}(t):=e^{\int_0^t(c+\gamma)(X_s)ds+N^{\gamma,\hat{b}^{H}}_tds},\
\ t\ge0,
$$ and
$$
{\cal J}^l_t:=e_{\gamma}(t)v(X_t)1_{\{\tau_{E_l\cap D}>t\}}.
$$
By \cite[Theorem 1.15]{MSW} and \cite[Lemma A.6, Theorem A.8 and
Lemma A.9]{MMS}, similar to \cite[Corollary 5.2.1]{Fuku} we can
show that there exists a subsequence $\{k_l\}$ such that for q.e.
$x\in\mathbb{R}^d$,
\begin{equation}\label{may13}
P_x\left\{\lim_{l\rightarrow\infty}N^{\gamma,\hat{b}^H_{k_l}}_t=N^{\gamma,\hat{b}^H}_t\
{\rm uniformly\ on\ any\ finite\ interval\ of}\ t\right\}=1.
\end{equation}
For simplicity, we still use $\{k\}$ to denote the subsequence
$\{k_l\}$.
 Letting $k\rightarrow\infty$, we obtain by (\ref{12}), (\ref{may2})--(\ref{may4}), (\ref{Oshi})--(\ref{may13}) and
\cite[Theorem 2.8 and (22)]{CMS} that there exists a sequence of
stopping times $\{\varpi_n\uparrow\infty\}$ such that $({\cal
J}^l_{t\wedge\varpi_n})_{t\ge 0}$ is a martingale under
$P^{\gamma}_x$  for $dx$-a.e. $x\in E_l\cap D$ and $n\in
\mathbb{N}$. Then, for $dx$-a.e. $x\in E_l\cap D$ and $n\in
\mathbb{N}$, we have
\begin{equation}\label{pp1}
v(x)=E_x^{\gamma}[e_{\gamma}(t\wedge\varpi_n)v(X_{t\wedge\varpi_n})1_{\{\tau_{E_l\cap
D}>t\wedge\varpi_n\}}].
\end{equation}
Letting $t,n\rightarrow\infty$, we obtain by (\ref{14}),
(\ref{pp1}) and \cite[Lemma 2.3(1)]{S} that $v(x)=0$ for $dx$-a.e. $x\in E_l\cap D$. Since
$l\in \mathbb{N}$ is arbitrary and $v|_D\in
 C(D)$, we obtain that $v\equiv0$ on $\mathbb{R}^d$. The proof is complete.\hfill\fbox

\bigskip

{ \noindent {\bf\large Acknowledgments} \quad   This work was
supported by Natural Sciences and Engineering Research
Council of Canada. }
%We thank the referee for the careful reading of our paper and all of the insightful
%comments that greatly improved the presentation of the paper.

%\begin{rem}
%From the above proofs and \cite[Theorem 1.3]{CH}, we can see that
%the positive constant $M$ of Theorem \ref{thm1} depends only on
%$d,\alpha,a,\|b\|_{L^{2p}},\|\hat{b}\|_{L^{2p}}$ and $D$.

%\end{rem}

\end{document}